\documentclass[12pt]{amsart}
\usepackage[latin1]{inputenc}
\usepackage{amsfonts,amsxtra,amssymb,amsthm,amsmath,amscd,mathrsfs, epsfig,url}
\def \Z {{\mathbb Z}}

\def\leq{\leqslant}

\def\geq{\geqslant}

\renewcommand{\Re}{{\rm Re\,}}

\theoremstyle{plain}

\theoremstyle{remark}

\theoremstyle{definition}

\numberwithin{equation}{section}

\begin{document}


\title{\bf The average size of Ramanujan sums over quadratic  number fields(II) }

\author{ Wenguang Zhai}




\thanks{This work is  supported by the National Natural Science Foundation of China(Grant No.  11971476).}

\maketitle

\addtocounter{footnote}{1}

Abstract:   In this paper we study Ramanujan sums $c_{\bf m}(\bf n)$, where
 $ {\bf m}$ and ${\bf n}$ are integral ideals in an arbitrary quadratic number field.  
We give some new results about the asymptotic behavior of sums of $c_{\bf m}(\bf n)$  over both $ {\bf m}$ and
$ {\bf n}$.

\section{\bf Introduction and statements of results}\
\subsection{\bf Ramanujan sums over rationals}\
For any positive integers $m$ and $n$, the classical Ramanujan sum $c_m(n)$ is defined by (see, for example, Kr\"{a}tzel \cite{Kr})
\begin{equation}
c_m(n):=\sum_{\stackrel{1\leq j\leq m}{gcd(j,m)=1}}e\left(\frac{jn}{m}\right)=\sum_{d|gcd(m,n)}d\mu\left(\frac md\right),
\end{equation}
 where $e(z):=e^{2\pi iz}$ and $\mu(\cdot)$ is the M\"{o}bius function. The Ramanujan sum is an interesting and important object in number theory and there
 are lots of papers in this area.

 In 2012,  Chan and Kumchev \cite{CK} dealt
with the question of the average order of $c_m(n)$ with respect to both variables
$m$ and $n$. Let $Y\geq X\geq 3$ be two large real numbers and $k\geq 1$ be a fixed integer. Define
\begin{eqnarray}
C_k(X, Y):&&=\sum_{1\leq n\leq Y}\left(\sum_{1\leq m\leq X} c_m(n)\right)^k.
\end{eqnarray}

When $k=1$,
they proved that the asymptotic formula
\begin{eqnarray}
C_1(X, Y)=Y-\frac{3}{2\pi^2}X^2+O(XY^{1/3}\log X)+O(X^3Y^{-1})
\end{eqnarray}
holds, which implies that if $Y\asymp X^{ \delta}$ then
\begin{eqnarray}
S_1(X, Y) \sim\left\{\begin{array}{ll}
Y,&\mbox{if  $ \delta>2,$}\\
 -\frac{3}{2\pi^2}X^2,&\mbox{if $1< \delta<2.$}
\end{array}\right.
\end{eqnarray}

When $k=2,$ Chan and Kumchev proved that if $Y\geq x^2(\log B)^B $ for some fixed $B>0$, then
\begin{equation}
C_2(X,Y)=\frac{YX^2}{2\zeta(2)}+O(X^4+XY\log X).
\end{equation}

\subsection{\bf Ramanujan sums over quadratic number fields}\
Suppose ${\Bbb F}/{\Bbb Q}$ is a number field of degree ${\bf d}\geq 2$ and its ring of algebraic integers is denoted by $\mathcal{O}_{\Bbb F}.$
For any nonzero integral ideal $\mathcal{I}$ in $\mathcal{O}_{\Bbb F}$, the M\"{o}bius function is defined as follows(see, for example \cite{He}, Page 100):
$\mu(\mathcal{I})=0$ if there exists a prime ideal $\mathcal{P}$ such that $\mathcal{P}^2$ divides $\mathcal{I},$ and
$\mu(\mathcal{I})=(-1)^r$ if $\mathcal{I}$ is a product of $r$ distinct prime ideals. For any ideal $\mathcal{I},$ the norm
of $\mathcal{I}$ is denoted by $N(\mathcal{I}).$

 For two nonzero integral ideals ${\bf m}$ and ${\bf n},$ the Ramanujan sum
is defined by
\begin{equation}
c_{\bf m}(\bf n):=\sum_{\stackrel{  {\bf d}\in \mathcal{O}_{\Bbb F}}{  {\bf d}| {\bf m},  {\bf d}| {\bf n}     }}
N( {\bf d})\mu\left(\frac{ {\bf m}}{ {\bf d}}\right),
\end{equation}
which is an analogue of (1.1). The definition (1.6)  of Ramanujan sums can be considered in the much
more general context of arbitrary arithmetical semigroups. See, for example, Grytczuk
\cite{Gr}  and the monograph by Knopfmacher \cite{Kn}.

For each $n\geq 1$, let $a_{\Bbb F}(n)$ denote the number of integral ideals  $\mathcal{I}$ in $\mathcal{O}_{\Bbb F}$ such that $N(\mathcal{I})=n.$ Then
we have
\begin{equation}
 \sum_{n\leq x}a_{\Bbb F}(n)=\rho_{\Bbb F}x+P_{\Bbb F}(x), \ \ P_{\Bbb F}(x)= O(x^{\frac{{\bf d}-1}{{\bf d}+1}}),
\end{equation}
where $\rho_{\Bbb F}$ is a constant depending only on ${\Bbb F}.$ The asymptotic formula (1.7) is a classical result of  Landau  (see \cite{Lan}),
which was improved by many authors(see, for example, M\"{u}ller \cite{Mu}, Nowak \cite{No3}).

Let $X\geq 3$ and $Y\geq 3$ be two large real numbers and $k\geq 1$ be a fixed integer. Define
\begin{eqnarray}
C_{{\Bbb F},k}(X, Y):=\sum_{1\leq N({\bf n})\leq Y}\left(\sum_{1\leq N({\bf m})\leq X}c_{\bf m}({\bf n})\right)^k,
\end{eqnarray}
which is an analogue of the sum $C_k(X, Y)$ defined by (1.2).

In \cite{No1}, W. G.  Nowak proved that if ${\Bbb F}$ is a  fixed quadratic number field, then the asymptotic formula
\begin{equation}
C_{{\Bbb F}, 1}(X, Y)\sim \rho_{\Bbb F}Y
\end{equation}
holds provided that $Y>X^\delta$ for some
$ \delta>  1973/820=2.40609\cdots. $
In \cite{No2}, W. G. Nowak had considered the Gaussian field  case ${\Bbb F}={\Bbb Q}(i)$, where he proved that the asymptotic formula (1.9) holds
 provided that  $Y>X^\delta$ for $ \delta>29/12=2.41\dot{6}.$

In \cite{Zh}, the author proved that the asymptotic formula (1.9)   holds
 provided that  $Y>X^\delta$ for $ \delta> 79/34=2.3235\cdots.$  It was also proved that (1.9) holds on average  for $2<\delta\leq 79/34.$

In this paper, we shall  prove that (1.9) holds  provided that  $Y>X^\delta$ for any $ \delta> 2.$ More precisely, we have the following Theorem 1.

{\bf Theorem 1.}  {\it Let ${\Bbb F}$ be a fixed quadratic number field and $3\leq X<Y$ be two large real numbers. Then
  we have
\begin{eqnarray}
C_{{\Bbb F}, 1}(X, Y)=\rho_{\Bbb F}Y+
 O\left(XY^{1/2}(\log Y)^7 +X^{2 }\right).
\end{eqnarray}
}

{\bf Remark 1.} Let $\lambda=\lambda(t)$ be an increasing function such that
 $\lim_{t\rightarrow \infty}\lambda(t)=\infty$ and $\lambda(t)=o(\log t)$ as  $t\rightarrow\infty.$  If $Y\geq X^2(\log X)^{14}\lambda(X),$ then   (1.9) holds.

We can also study the sum for  $k=2.$ In this case we have the following Theorem 2, which is an analogue and generalization of (1.5).

{\bf Theorem 2.}  {\it Let ${\Bbb F}$ be a fixed quadratic number field and $3\leq X<Y$ be two large real numbers. If $Y>X^2,$ then
\begin{eqnarray}
C_{\Bbb F,2}(X,Y)
&&=\frac{\rho^2_{\Bbb F}}{2\zeta_{\Bbb F}(2)}X^2Y+\frac{\zeta_{\Bbb B}(0)\rho^2_{\Bbb F}}{4\zeta^2_{\Bbb F}(2)}X^4 \\
 &&\ \ \ +O\left(X^{\frac{24}{5}}Y^{-\frac 25}+X^{2}Y^{\frac 23} \log^5 Y +X^{\frac{3}{2}}Y \log^3 Y
\right)\nonumber
\end{eqnarray} }

{\bf Remark 2.}  In \cite{Zh}, the method of exponential sums was applied. However, in this paper, we apply the method of complex integration.

{\bf Notation.} Throughout this paper, we use the following notations. ${\Bbb N}$,  ${\Bbb Q}$  and ${\Bbb F}$  denote the set of positive integers, the set of rational numbers and  a number field of degree ${\bf d}\geq 2,$ respectively. We say
$n$ is a half integer if $n-1/2\in {\Bbb N}.$
For each $n\in {\Bbb N},$ $a_{\Bbb F}(n)$ denotes the number of integral ideals $\mathcal{I}$
such that $N(\mathcal{I})=n,$ $\tau_{\ell }(n)$ denotes the number of ways $n$ can be written as a product of $\ell$ positive integers, $\tau(n)=\tau_2(n)$ is the well-known Dirichlet divisor function. $\zeta(s)$ is the Riemann zeta-function, and $\zeta_{\Bbb F}(s)$
is the Dedekind zeta-function of the field ${\Bbb F}.$ $\varepsilon$ always denotes a small positive constant, which may be different at different places.

\section{\bf Preliminary Lemmas}

We   suppose that ${\Bbb F}$  is a fixed number field of   degree $d\geq 2.$
The Dedekind zeta function of ${\Bbb F}$ is defined by
\begin{equation}
\zeta_{\Bbb F}(s):=\sum_{\stackrel{\mathcal{I}\in \mathcal{O}_{\Bbb F}}{\mathcal{I}\not= 0}}\frac{1}{N^s(\mathcal{I})}\ \ (\Re (s)>1).
\end{equation}
Then
\begin{eqnarray}
\zeta_{\Bbb F}(s)=\sum_{n=1}^\infty \frac{a_{\Bbb F}(n)}{n^s}\ \ (\Re(s)>1),
\end{eqnarray}
where $a_{\Bbb F}(n)$ denotes the number of integral ideals $\mathcal{I}$
such that $N(\mathcal{I})=n.$

We have
$$\frac{1}{\zeta_{\Bbb F}(s)}= \sum_{\stackrel{ {\bf n}\in \mathcal{O}_{\Bbb F}}{{\bf n}\not= 0}} \frac{\mu({\bf n})}{N^s({\bf n})}\ \ (\Re(s)>1).$$

{\bf Lemma 2.1.}   {\it Suppose ${\Bbb F}$  is a fixed number field of   degree $d\geq 2.$ Then we have the functional equation
$$\zeta_{\Bbb F}(s)=\chi_{\Bbb F}(s)\zeta_{\Bbb F}(1-s) $$
such that the estimate
$$\chi_{\Bbb F}(s)\ll (|t|+1)^{d(\frac 12-\sigma)}\ \ (|t|\rightarrow \infty)$$
holds in any fixed critical strip.
}

\begin{proof}
See, for example, Iwaniec and Kowalski \cite {IK}.
\end{proof}

{\bf Lemma 2.2.} {\it Suppose ${\Bbb F}$  is a fixed number field of   degree $d\geq 2.$

If $\sigma\geq 1,$ then
\begin{equation}
\zeta_{\Bbb F}(\sigma+it)\ll_{\Bbb F} \log(|t|+2).
\end{equation}

If $\sigma\geq 1,$ then
\begin{equation}
\frac{1}{\zeta_{\Bbb F}(\sigma+it)}\ll_{\Bbb F} \log(|t|+2).
\end{equation}

If $0\leq \sigma\leq 1,$ then
\begin{equation}
\zeta_{\Bbb F}(\sigma+it)\ll_{\Bbb F} (|t|+2)^{\frac{d(1-\sigma)}{2}}\log(|t|+2).
\end{equation}

If $-2\leq \sigma\leq 0,$ then
\begin{equation}
\zeta_{\Bbb F}(\sigma+it)\ll_{\Bbb F} (|t|+2)^{d(\frac 12-\sigma)}\log(|t|+2).
\end{equation}
}

\begin{proof}
Let $x=(|t|+2)^{d+1}.$ Suppose $\sigma>1.$ We write
\begin{equation}
\zeta_{\Bbb F}(\sigma+it)=\Sigma_1+\Sigma_2,
\end{equation}
where
$$\Sigma_{1}:=\sum_{n\leq x}\frac{a_{\Bbb F}(n)}{n^{\sigma+it}},
\ \ \ \Sigma_{2}:=\sum_{n> x}\frac{a_{\Bbb F}(n)}{n^{\sigma+it}}.
$$

Let $A_{\Bbb F}(u)=\sum_{n\leq u}a_{\Bbb F}(n).$
By   partial summation we have
\begin{eqnarray}
\Sigma_2=&&\int_x^\infty \frac{dA_{\Bbb F}(u)}{u^{\sigma+it}}
=\rho_{\Bbb F}\int_x^\infty \frac{du}{u^{\sigma+it}}+
\int_x^\infty \frac{dP_{\Bbb F}(u)}{u^{\sigma+it}}\\
&&=\rho_{\Bbb F}\frac{x^{1-\sigma-it}}{\sigma-1+it}
-\frac{P_{\Bbb F}(x)}{x^{\sigma+it}}
+(\sigma+it)
\int_x^\infty \frac{P_{\Bbb F}(u)}{u^{\sigma+1+it}}du.\nonumber
\end{eqnarray}

By the estimate $P_{\Bbb F}(u)\ll u^{(d-1)/(d+1)}$ we see that
(2.8) is valid for $\sigma>(d-1)/(d+1).$ So for $\sigma\geq 1$ we have
(recalling the definition of $x$)
\begin{eqnarray}
\Sigma_2\ll 1+ (|t|+2) x^{-\frac{2}{d+1}}\ll 1.
\end{eqnarray}

For $\Sigma_1$ we have by partial summation that
\begin{equation}
\Sigma_1\ll \sum_{n\leq x}\frac{a_{\Bbb F}(n)}{n }\ll \log x\ll \log(|t|+2)\ (\sigma\geq 1).
\end{equation}

Now the estimate (2.3) follows from (2.7)-(2.10). The proof of (2.4) is similar and easier.

By (2.3) and Lemma 2.1 we get
\begin{equation}
\zeta_{\Bbb F}(it)\ll_{\Bbb F} (|t|+2)^{\frac{d}{2}}\log(|t|+2).
\end{equation}
So (2.5) follows from (2.11) and
$\zeta_{\Bbb F}(1+it)\ll_{\Bbb F} \log(|t|+2)$ with the help
of Phragmen-Lindel\"{o}f principle.
The estimate (2.6) follows from (2.3) and Lemma 2.1.
\end{proof}

{\bf Lemma 2.3.} {\it  Suppose ${\Bbb F}$ is a quadratic number field. Then the estimate
$$\int_{-U}^U |\zeta_{\Bbb F}(\sigma+it)|^2dt\ll U(\log U)^4\ \ (U\geq 2)$$
holds uniformly for $1/2\leq \sigma\leq 1.$
}

\begin{proof}
Suppose  ${\Bbb F}$ is a quadratic number field, then there  exists a real primitive Dirichlet Character $\chi_D$ of modulo $|D|$  such that
$ \zeta_{\Bbb F}(s)=\zeta(s)L(s, \chi_D),$
where $L(s, \chi_D)$ is the Dirichlet L-function corresponding to $\chi_D.$
Now Lemma 2.3 follows from the fourth power moment of $\zeta(s)$ and
the fourth power moment of $ L(s, \chi_D)$.  \end{proof}

{\bf Lemma 2.4.} {\it  Suppose ${\Bbb F}$ is a quadratic number field. Then the estimate
$$\int_{-U}^U |\zeta_{\Bbb F}(\sigma+it)|^4dt\ll U \ \ (U\geq 2)$$
holds uniformly for $2/3\leq \sigma\leq 1.$
}

\begin{proof}
We have
\begin{eqnarray*}
\int_1^U \left|\zeta(\frac 23+it)\right|^8dt\ll U \ \ (U\geq 2).
\end{eqnarray*}
In \cite{Iv} we can find even much stronger results. Similarly, we have
\begin{eqnarray*}
\int_1^U \left|L(\frac 23+it,\chi_D)\right|^8dt\ll U \ \ (U\geq 2).
\end{eqnarray*}
So Lemma 2.4 follows from the above two estimates and Cauchy's inequality.
 \end{proof}

{\bf Lemma 2.5.} {\it  Suppose ${\Bbb F}$ is a quadratic number field. Then for
$1/2\leq \sigma\leq 1$
we have the estimate
$$\zeta_{\Bbb F}(\sigma+it)\ll_{\Bbb F}
(|t|+2)^{\frac 23 (1-\sigma)}\log(|t|+2). $$
}

\begin{proof}
 This estimate follows from (2.3) with $d=2$, the well-known bound
 $$\zeta_{\Bbb F}( 1/2+it)\ll_{\Bbb F}
(|t|+2)^{\frac 13 } $$
and the  Phragmen-Lindel\"{o}f principle.
 \end{proof}

{\bf Lemma 2.6. } {\it Let $k\geq 2$ be a fixed integer, and  $f(n_1,\cdots, n_k)$ is a multivariable  arithmetic function
such that its Dirichlet series
\begin{eqnarray}
F(s_1,\cdots, s_k)=\sum_{n_1= 1}^\infty\cdots \sum_{n_k= 1}^\infty\frac{f(n_1,\cdots, n_k)}{n_1^{s_1}\cdots  n_k^{s_k}}
\end{eqnarray}
is absolutely convergent for $\Re (s_j)>\sigma_j
\ (j=1,\cdots, k),$ where   $\sigma_1>0,\cdots, \sigma_k>0.$
Suppose $x_1,\cdots, x_k, T_1,\cdots,  T_k\geq 5$  are  parameters such that $x_j\notin {\Bbb N}\ (j=1,\cdots, k)$,   and define
$$b_j =\sigma_j+\frac{1}{\log x_j}, \ \  \ (j=1,\cdots, k).$$
  Then we have
\begin{eqnarray}
&&\ \ \ \ \sum_{n_1\leq x_1}\cdots \sum_{n_k\leq x_k}f(n_1,\cdots, n_k)\\  &&
 =\frac{1}{(2\pi i)^k}\int_{b_1-iT_1}^{b_1+iT_1} \cdots \int_{b_k-iT_k}^{b_k+iT_k}\frac{F(s_1,\cdots, s_k)x_1^{s_1}\cdots x_k^{s_k}}{s_1\cdots s_k} \mbox{d}s_k \cdots \mbox{d}s_1 \nonumber  \\
&& \ \ \ \ \ +O(x_1^{\sigma_1}\cdots x_k^{\sigma_k}E),\nonumber
\end{eqnarray}
where
 \begin{equation}
E :=\sum_{j=1}^k\sum_{n_1=1}^\infty\cdots \sum_{ n_k= 1}^\infty \frac{|f(n_1,\cdots, n_k)|}{ n_1^{b_1}\cdots n_k^{b_k}  }
\times\frac{1}{T_j|\log\frac{x_j}{n_j}|+1}.
\end{equation}
}

 \begin{proof}

For $b>0, a>0, T>1,$ we have
\begin{equation}
\frac{1}{2\pi i}\int_{b-iT}^{b+iT}\frac{a^{s}}{s}\mbox{d}s=\delta(a)
+O\left( \frac{a^b }{T|\log a|+1} \right),
\end{equation}
where $\delta(a)=1$ when $a>1,$ and $\delta(a)=0$ when $0<a<1.$ See for example, Chapter 2 in the second part of Tenenbaum \cite{Te}.

The Dirichlet series (2.12) is absolutely   convergent for $\Re(s_j)>\sigma_j\ (j=1,\cdots, k).$ So we  have by (2.15) that
\begin{eqnarray}
&&\ \ \ \ \ \ \ \ \frac{1}{(2\pi i)^k}\int_{b_1-iT_1}^{b_1+iT_1} \cdots \int_{b_k-iT_k}^{b_k+iT_k}F(s_1,\cdots, s_k)\frac{x_1^{s_1}\cdots x_k^{s_k}}{s_1\cdots s_k}   \mbox{d}s_k\cdots  \mbox{d}s_1  \\
&&=\frac{1}{(2\pi i)^k}\int_{b_1-iT_1}^{b_1+iT_1}\cdots  \int_{b_k-iT_k}^{b_k+iT_k}
\sum_{n_1= 1}^\infty\cdots \sum_{n_k= 1}^\infty\frac{f(n_1,\cdots, n_k)}{n_1^{s_1}\cdots  n_k^{s_k}}
 \frac{x_1^{s_1}\cdots x_k^{s_k}}{s_1\cdots s_2}   \mbox{d}s_k\cdots  \mbox{d}s_1  \nonumber\\
&&=\sum_{n_1= 1}^\infty\cdots \sum_{n_k= 1}^\infty
\frac{f(n_1,\cdots, n_k) }{(2\pi i)^k}
\int_{b_1-iT_1}^{b_1+iT_1}\cdots  \int_{b_k-iT_k}^{b_k+iT_k}
\frac{1}{n_1^{s_1}\cdots  n_k^{s_k}}
 \frac{x_1^{s_1}\cdots x_k^{s_k}}{s_1\cdots s_k}   \mbox{d}s_k\cdots  \mbox{d}s_1  \nonumber\\
&&=\sum_{n_1= 1}^\infty\cdots \sum_{n_k= 1}^\infty f(n_1,\cdots, n_k)
 \prod_{j=1}^k\left(\frac{1}{2\pi i}\int_{b_j-iT_j}^{b_j+iT_j}\left(\frac{x_j}{n_j}\right)^{s_j} \frac{\mbox{d}s_j}{s_j}\right)\nonumber\\
&&=\sum_{n_1= 1}^\infty\cdots \sum_{n_k= 1}^\infty f(n_1,\cdots, n_k)
 \prod_{j=1}^k\left(\delta\left(\frac{x_j}{n_j}\right)
 +E_j\left(\frac{x_j}{n_j}\right)\right), \nonumber
\end{eqnarray}
where
$$E_j\left(\frac{x_j}{n_j}\right)=O\left(\frac{(x_j/n_j)^{b_j}}{T_j\left| \log\frac{x_j}{n_j}\right|+1}\right)\ (j=1,\cdots,k).$$

For any integer $n_j>0$ we have
\begin{eqnarray*}
&&\delta\left(\frac {x_j}{n_j}\right)\leq \left(\frac{x_j}{n_j}\right)^{b_j}\ll \frac{x_j^{\sigma_j}}{n_j^{b_j}},\\
&&E_j\left(\frac{x_j}{n_j}\right)\ll \frac{x_j^{\sigma_j}}{n_j^{b_j}}\times
\frac{1}{T_j\left| \log\frac{x_j}{n_j}\right|+1}\ (j=1,\cdots,k).
\end{eqnarray*}

Thus we have
\begin{eqnarray}
&&\prod_{j=1}^k\left(\delta\left(\frac{x_j}{n_j}\right)
 +E_j\left(\frac{x_j}{n_j}\right)\right)\\
&&=\delta\left(\frac{x_1}{n_1}\right)\cdots \delta\left(\frac{x_k}{n_k}\right)
+\sum_{j=1}^kO\left(\frac{x_1^{\sigma_1}\cdots x_k^{\sigma_k}}
{n_1^{b_1}\cdots n_k^{b_k}} \times
\frac{1}{T_j\left| \log\frac{x_j}{n_j}\right|+1} \right).\nonumber
\end{eqnarray}

Now Lemma 2.4 follows from (2.16) and (2.17) by noting that
\begin{eqnarray*}
 \sum_{n_1= 1}^\infty\cdots \sum_{n_k= 1}^\infty f(n_1,\cdots, n_k)\delta\left(\frac{x_1}{n_1}\right)\cdots \delta\left(\frac{x_k}{n_k}\right)=\sum_{n_1\leq x_1}\cdots  \sum_{n_k\leq x_k} f(n_1,\cdots, n_k).
\end{eqnarray*}
\end{proof}

\section{\bf Some special Dirichlet series}

\subsection{\bf Dirichlet series involving $c_{{\bf m}}({\bf n}) $}\
For fixed $k\geq 1,$ we define a multivariate arithmetic function
$f({\bf m}_1,\cdots, {\bf m}_k, {\bf n})$
over the number field  ${\Bbb F}$ by
\begin{equation}
f({\bf m}_1,\cdots, {\bf m}_k, {\bf n}):=c_{{\bf m}_1}({\bf n})\cdots c_{{\bf m}_k}(\bf n).
\end{equation}
Note that when $k=1,$ $f({\bf m}_1,{\bf n})=c_{{\bf m}_1}(\bf n).$

Suppose $s_1, \cdots, s_k, w\in {\Bbb C}$ with $\Re(s_j)>1(j=1,\cdots,k), \Re (w)>2.$ Define the Derichlet series
\begin{equation}
\mathcal{F}(s_1,\cdots,s_k,w):=\sum_{{\bf m}_1, \cdots, {\bf m}_k,{\bf n}}
\frac{f({\bf m}_1,\cdots, {\bf m}_k, {\bf n})}{ N^{s_1}({\bf m}_1)\cdots  N^{s_k}({\bf m}_k) N^w({\bf n})}.
\end{equation}

For any $\theta\in {\Bbb C}$ and  any non-zero integral ideal ${\bf n}$, we define the weighted divisor function
\begin{equation}
\sigma_{\theta}({\bf n}):= \sum_{{\bf d}|{\bf n}} N^{\theta}(\bf d).
\end{equation}

We have the following Lemma 3.1.

{\bf Lemma 3.1.}  {\it Suppose $\theta_1,  \theta_2, w\in {\Bbb C}.$

If     $\Re{w}>\max(1, 1+\Re(\theta_1)), $ then we have the  identity
\begin{equation}
\sum_{{\bf n}}\frac{\sigma_{\theta_1}({\bf n})}{N^w({\bf n})}=\zeta_{\Bbb F}(w)\zeta_{\Bbb F}(w-\theta_1).
\end{equation}

If     $\Re{w}>\max(1, 1+\Re(\theta_1), 1+\Re(\theta_1),  1+\Re(\theta_1+\theta_2)   ), $ then we have   the Ramanujan's identity
\begin{equation}
\sum_{{\bf n}}\frac{\sigma_{\theta_1}({\bf n}) \sigma_{\theta_2}({\bf n})    }{N^w({\bf n})}=
\frac{\zeta_{\Bbb F}(w)\zeta_{\Bbb F}(w-\theta_1)\zeta_{\Bbb F}(w-\theta_2) \zeta_{\Bbb F}(w-\theta_1-\theta_2) }
 {\zeta_{\Bbb F}(2w-\theta_1-\theta_2)}.
\end{equation}

}

\begin{proof}
The formula (3.4) follows from the definitions of $\sigma_\theta(\cdot)$ and the Dedekind zeta-function (2.1). The formula (3.5) can be proved in the same
way as the proof of the formula (1.3.3) in \cite{Ti}. We omit the details.
\end{proof}

For the function $\mathcal{F}(s_1,\cdots,s_k,w),$ we then have the following

{\bf Proposition 3.1.}   {\it   Suppose $w, s_1, s_2\in {\Bbb C}.$

If $ \Re(w)>1, \Re(w+s_1)>2, $ then
\begin{eqnarray}
\mathcal{F}(s_1,w)=\frac{\zeta_{\Bbb F}(w)\zeta_{\Bbb F}(w+s_1-1)}{\zeta_{\Bbb F}(s_1)},
\end{eqnarray}

If $\Re(w)>1, \Re(w+s_1)>2,  \Re(w+s_2)>2, \Re(w+s_1+s_2)>3,   $ then
\begin{eqnarray}
\ \ \ \mathcal{F}(s_1,s_2,w)=\frac{ \zeta_{\Bbb F}(w)   \zeta_{\Bbb F}(w+s_1-1)  \zeta_{\Bbb F}(w+s_2-1)  \zeta_{\Bbb F}(w+s_1+s_2-2)    }
{\zeta_{\Bbb F}(s_1)\zeta_{\Bbb F}(s_2)  \zeta_{\Bbb F}(2w+s_1+s_2-2)   }.
\end{eqnarray}

}

\begin{proof}
Obviously we can rewrite the formula (3.2) in the form
\begin{equation}
\mathcal{F}(s_1,\cdots,s_k,w)=\sum_{\bf n}\frac{1}{N^w({\bf n})} \prod_{j=1}^k \left( \sum_{{\bf m}_j}\frac{c_{{\bf m}_j}({\bf n})}{N^{s_j}({\bf m}_j)}   \right).
\end{equation}

Suppose $s\in {\Bbb C}$ such that $\Re(s)>1.$   By (1.6)   we have
\begin{eqnarray}
\sum_{{\bf m}}\frac{c_{{\bf m}}({\bf n})}{N^{s}({\bf m})}&& =\sum_{{\bf m}}\frac{1}{N^{s}({\bf m})}\sum_{{\bf d}|{\bf m}, {\bf d}|{\bf n}}N({\bf d})\mu\left(\frac{{\bf m}}{{\bf d}}\right)\\
&&=\sum_{{\bf d}|{\bf n}}N^{1-s}({\bf d})\sum_{{\bf m}^{*}}\frac{\mu({\bf m}^{*})}{N^s({\bf m}^{*})}\nonumber\\
&&=\frac{\sigma_{1-s}({\bf n})}{\zeta_{\Bbb F}(s)}.\nonumber
\end{eqnarray}

From (3.8) and (3.9) we get
\begin{equation}
\mathcal{F}(s_1,\cdots,s_k,w)=\frac{1}{\zeta_{\Bbb F}(s_1)\cdots \zeta_{\Bbb F}(s_k)}
\sum_{\bf n}\frac{ \sigma_{1-s_1}({\bf n}) \cdots   \sigma_{1-s_k}({\bf n})     }{N^w({\bf n})}.
\end{equation}

Now Proposition 3.1 follows from  (3.10) and Lemma 3.1.
\end{proof}

\subsection{\bf Dirichlet series involving
$c_{\bf m}^{*}(\bf n) $}\
For non-zero integral ideals ${\bf m}$ and ${\bf n},$ define
\begin{equation}
c_{\bf m}^{*}(\bf n):=\sum_{\stackrel{  {\bf d}\in \mathcal{O}_{\Bbb F}}{  {\bf d}| {\bf m},  {\bf d}| {\bf n}     }}
N( {\bf d})\left|\mu\left(\frac{ {\bf m}}{ {\bf d}}\right)\right|.
\end{equation}

It is easily seen that
\begin{equation}
 |c_{\bf m}(\bf n)|\leq c_{\bf m}^{*}(\bf n).
 \end{equation}

Suppose $s\in {\Bbb C}$ such that $\Re(s)>1$. It is easy to see that
\begin{eqnarray}
\sum_{{\bf m}_{1}}\frac{|\mu({\bf m}_{1})|}{N^s({\bf m}_{1})}=
\frac{\zeta_{\Bbb F}(s)}{\zeta_{\Bbb F}(2s)}.
\end{eqnarray}
So for any non-zero integral ideal
${\bf n}$, we have for $\Re(s)>1$ that
\begin{eqnarray}
\sum_{{\bf m}}\frac{c_{{\bf m}}^{*}({\bf n})}{N^{s}({\bf m})}&& =\sum_{{\bf m}}\frac{1}{N^{s}({\bf m})}\sum_{{\bf d}|{\bf m}, {\bf d}|{\bf n}}N({\bf d})\left|\mu\left(\frac{{\bf m}}{{\bf d}}\right)\right|\\
&&=\sum_{{\bf d}|{\bf n}}N^{1-s}({\bf d})\sum_{{\bf m}_{1}}\frac{|\mu({\bf m}_{1})|}{N^s({\bf m}_{1})}\nonumber\\
&&=\frac{\zeta_{\Bbb F}(s)}{\zeta_{\Bbb F}(2s)}\sigma_{1-s}({\bf n}).\nonumber
\end{eqnarray}

\section{\bf Estimates of some sums}

Suppose $X, Y, T\geq 3$  are large  real numbers such that both $X$ and $Y$ are
half integers and $X\leq Y.$
  Let $\sigma_0=1+1/\log X$ and $k\geq 1$ be a fixed integer. In this section we shall estimate the  sums
 $E_{j,k}(X,T)\ (j=1,2,\cdots,k)$ and $\mathfrak{E}_{k}(Y,T),$  which are defined by
\begin{eqnarray*}
&& E_{j,k}(X,T):=
 \sum_{\stackrel{{\bf m}_1,\cdots, {\bf m}_k\in \mathcal{O}_{\Bbb F} }  {{\bf n}\in \mathcal{O}_{\Bbb F}}}
\frac{c_{{\bf m}_1}^{*} ({\bf n}) \cdots  c_{{\bf m}_k}^{*} ({\bf n})      }{N^{\sigma_0}({\bf m}_1) \cdots N^{\sigma_0}({\bf m}_k)    N^{\sigma_0}({\bf n})     }
\times  \frac{1}{T\left|\log \frac{X}{N({\bf m}_j)}\right|+1}
 \end{eqnarray*}
 and
 \begin{eqnarray*}
 &&\mathfrak{E}_{k}(Y,T):=  \sum_{\stackrel{{\bf m}_1,\cdots, {\bf m}_k\in \mathcal{O}_{\Bbb F} }  {{\bf n}\in \mathcal{O}_{\Bbb F}}}
\frac{c_{{\bf m}_1}^{*} ({\bf n}) \cdots  c_{{\bf m}_k}^{*} ({\bf n})      }{N^{\sigma_0}({\bf m}_1) \cdots N^{\sigma_0}({\bf m}_k)    N^{\sigma_0}({\bf n})     }
 \times  \frac{1}{T\left|\log \frac{Y}{N({\bf n})}\right|+1}
\end{eqnarray*}
respectively. It is easy to see that
$E_{1,k}(X,T)=E_{2,k}(X,T)=\cdots =E_{k,k}(X,T).$
Thus it suffices to bound $E_{1,k}(X,T)$ and $\mathfrak{E}_{k}(Y,T)$.

\subsection{\bf An auxiliary estimate}

Suppose $T$ and $U$ are large real numbers such that $U$ is a half integer, $g(n)$ is a non-negative arithmetic function  such that
$g(n)\ll n^{\varepsilon}$ holds for any $\varepsilon>0.$
Define
$$G(s):=\sum_{n=1}^\infty\frac{g(n)}{n^s}\ \ (\Re(s)>1),$$
which is obviously absolutely convergent for $\Re(s)>1.$
Suppose $1<\sigma_1<11/10.$ Define
$$E(U,T;\sigma_1):=\sum_{n=1}^\infty \frac{g(n)}{n^{\sigma_1}}\times  \frac{1}{T|\log\frac Un|+1}.$$

We write
\begin{equation}
E(U,T;\sigma_1)=E_1(U,T;\sigma_1)+  E_2(U,T;\sigma_1) + E_3(U,T;\sigma_1),
\end{equation}
where
\begin{eqnarray*}
&&E_1(U,T;\sigma_1):=\sum_{n\leq U/2}  \frac{g(n)}{n^{\sigma_1}}\times  \frac{1}{T|\log\frac Un|+1},\\
&&E_2(U,T;\sigma_1):=\sum_{U/2<n\leq 2U}  \frac{g(n)}{n^{\sigma_1}}\times  \frac{1}{T|\log\frac Un|+1},\\
&&E_3(U,T;\sigma_1):=\sum_{n> 2U}  \frac{g(n)}{n^{\sigma_1}}\times  \frac{1}{T|\log\frac Un|+1}.
\end{eqnarray*}

Trivially we have
\begin{eqnarray}
&&E_1(U,T;\sigma_1) +E_3(U,T;\sigma_1)\ll \frac 1T\sum_{n=1}^\infty  \frac{g(n)}{n^{\sigma_1}}=\frac{G(\sigma_1)}{T}.
\end{eqnarray}

For $E_2(U,T;\sigma_1),$ we have
\begin{eqnarray}
E_2(U,T;\sigma_1)&&\ll U^{\varepsilon-1}\sum_{U/2<n\leq 2U}\frac{1}{T|\log\frac Un|+1} \ll U^\varepsilon T^{-1},
\end{eqnarray}
where   we used the estimate
\begin{eqnarray*}
 \sum_{U/2<n\leq 2U}\frac{1}{T|\log\frac Un|+1}\ll \frac{U\log U}{T},
\end{eqnarray*}
which is well-known in analytic number theory.

From (4.1)-(4.3) we get the estimate
\begin{equation}
E(U,T;\sigma_1)\ll \frac{G(\sigma_1)}{T}+\frac{U^\varepsilon}{T}.
\end{equation}

\subsection{\bf Estimate of $\mathfrak{E}_{k}(Y,T)$}

We write
\begin{eqnarray*}
 &&\mathfrak{E}_{k}(Y,T)=   \sum_{ {\bf n}\in \mathcal{O}_{\Bbb F}}\frac{1}{N^{\sigma_0}({\bf n})}\times \frac{1}{T\left|\log \frac{Y}{N({\bf n})}\right|+1}
\left( \sum_{  {\bf m}\in \mathcal{O}_{\Bbb F} }
\frac{c_{{\bf m}}^{*} ({\bf n})       }{N^{\sigma_0}({\bf m})      }\right)^k,
\end{eqnarray*}
which combining (3.14) gives
\begin{eqnarray}
 \mathfrak{E}_{k}(Y,T)&&=  \sum_{ {\bf n}\in \mathcal{O}_{\Bbb F}}\frac{1}{N^{\sigma_0}({\bf n})}\times \frac{1}{T\left|\log \frac{Y}{N({\bf n})}\right|+1}
\left(  \frac{\zeta_{\Bbb F}(\sigma_0)}{\zeta_{\Bbb F}(2\sigma_0)} \sigma_{1-\sigma_0}({\bf n})  \right)^k\\
&&=\frac{\zeta_{\Bbb F}^k(\sigma_0) }{\zeta_{\Bbb F}^k(2\sigma_0)} \sum_{ {\bf n}\in \mathcal{O}_{\Bbb F}}
\frac{ \sigma_{1-\sigma_0}^k({\bf n})  }{N^{\sigma_0}({\bf n})}\times \frac{1}{T\left|\log \frac{Y}{N({\bf n})}\right|+1}\nonumber\\
&&\leq \frac{\zeta_{\Bbb F}^k(\sigma_0) }{\zeta_{\Bbb F}^k(2\sigma_0)} \sum_{ {\bf n}\in \mathcal{O}_{\Bbb F}}
\frac{ \sigma_0^k({\bf n})  }{N^{\sigma_0}({\bf n})}\times \frac{1}{T\left|\log \frac{Y}{N({\bf n})}\right|+1}\nonumber
\end{eqnarray}
by noting that
$$\sigma_{1-\sigma_0}({\bf n})=
\sum_{{\bf d}|{\bf n}}(N({\bf d}))^{-\frac{1}{\log X}}
\leq \sum_{{\bf d}|{\bf n}}1 =\sigma_{ 0}({\bf n}).
$$

Suppose $s$ such that $\Re(s)>1$. Then we have
$$\sum_{\bf n}\frac{ \sigma_0({\bf n})}{N^s({\bf n})}=\zeta_{\Bbb F}^2(s)=
\sum_{n=1}^\infty\frac{a_{\Bbb F}\ast a_{\Bbb F}(n)}{n^s},$$
where
$$a_{\Bbb F}\ast a_{\Bbb F}(n)=\sum_{n=n_1n_2}a_{\Bbb F}(n_1) a_{\Bbb F}(n_2).$$

So for $\Re(s)>1$ we can write
\begin{eqnarray}
&&\ \ \ \ \ G_1(s):=\sum_{ {\bf n}\in \mathcal{O}_{\Bbb F}}
\frac{ \sigma_0^k({\bf n})  }{N^{s}({\bf n})}= \sum_{n=1}^\infty
\frac{g_1(n)}{n^s},
\end{eqnarray}
where
$$g_1(n):=(a_{\Bbb F}\ast a_{\Bbb F}(n))^k   a_{\Bbb F}(n).$$
By the well-known bound $a_{\Bbb F}(n)\ll n^\varepsilon$ we get
that $g_1(n)\ll n^\varepsilon.$ So from (4.4)-(4.6) we get
\begin{eqnarray}
 \mathfrak{E}_{k}(Y,T)\ll \frac{\zeta_{\Bbb F}^k(\sigma_0) }{\zeta_{\Bbb F}^k(2\sigma_0)}\times \frac{1}{T}\left( G_1(\sigma_0)+Y^\varepsilon\right).
\end{eqnarray}

By Euler's product we have
\begin{eqnarray*}
&&\ \ \ \ \ \sum_{ {\bf n}\in \mathcal{O}_{\Bbb F}}
\frac{ \sigma_0^k({\bf n})  }{N^{s}({\bf n})}=\prod_{{\bf p}\in \mathcal{O}_{\Bbb F}}\left(1+\sum_{\alpha=1}^\infty
\frac{\sigma_0^k({\bf p}^\alpha) }{N^{s}({\bf p}^\alpha)}   \right)\\
&&=\prod_{{\bf p}\in \mathcal{O}_{\Bbb F}} \left(1-\frac{1}{N^s({\bf p})}\right)^{-2^k}
\prod_{{\bf p}\in \mathcal{O}_{\Bbb F}}\left(1+\sum_{\alpha=1}^\infty
\frac{\sigma_0^k({\bf p}^\alpha) }{N^{s}({\bf p}^\alpha)}   \right)\left(1-\frac{1}{N^s({\bf p})}\right)^{2^k}\nonumber\\
&&=(\zeta_{\Bbb F}(s))^{2^k}H(s),\nonumber
\end{eqnarray*}
where
$$ H(s)=\prod_{{\bf p}\in \mathcal{O}_{\Bbb F}}\left(1+\sum_{\alpha=1}^\infty
\frac{\sigma_0^k({\bf p}^\alpha) }{N^{s}({\bf p}^\alpha)}   \right)\left(1-\frac{1}{N^s({\bf p})}\right)^{2^k}.$$
It is easy to see that $H(s)$ is absolutely convergent for $\Re(s)>1/2.$ So
We have
\begin{eqnarray}
 G_1(\sigma_0)\ll (\zeta_{\Bbb F}(\sigma_0))^{2^k}.
\end{eqnarray}

From (2.3), (4.7)  and (4.8) we get
\begin{eqnarray}
  \mathfrak{E}_{k}(Y,T)\ll  \frac{(\log X)^{2^k+k}}{T}+ \frac{Y^\varepsilon}{T}\ll
  \frac{Y^\varepsilon}{T}
\end{eqnarray}
by noting that $X\leq Y.$

\subsection{\bf Estimate of $E_{1,k}(X,T)$}

 By (3.14) and the definition of $c_{\bf m}^{*}({\bf n})$ in last section we have
\begin{eqnarray}
&&\ \ \ \ \ \ \  E_{1,k}(X,T)\\
&&=   \sum_{ {\bf n}\in \mathcal{O}_{\Bbb F}}\frac{1}{N^{\sigma_0}({\bf n})}\sum_{{\bf m}_1\in \mathcal{O}_{\Bbb F}}
\frac{1}{N^{\sigma_0}({\bf m}_1)} \frac{c_{{\bf m}_1}^{*}({\bf n})}{T\left|\log \frac{X}{N({\bf m}_1)}\right|+1}
\left( \sum_{  {\bf m}\in \mathcal{O}_{\Bbb F} }
\frac{c_{{\bf m}}^{*} ({\bf n})       }{N^{\sigma_0}({\bf m})      }\right)^{k-1}\nonumber\\
&&=\frac{\zeta_{\Bbb F}^{k-1}(\sigma_0) }{\zeta_{\Bbb F}^{k-1}(2\sigma_0)}\sum_{ {\bf n}\in \mathcal{O}_{\Bbb F}}\frac{\sigma_{1-\sigma_0}^{k-1}({\bf n})}{N^{\sigma_0}({\bf n})}\sum_{{\bf m}_1\in \mathcal{O}_{\Bbb F}}\frac{1}{N^{\sigma_0}({\bf m}_1)} \frac{c_{{\bf m}_1}^{*}({\bf n})}{T\left|\log \frac{X}{N({\bf m}_1)}\right|+1}\nonumber\\
&&=\frac{\zeta_{\Bbb F}^{k-1}(\sigma_0) }{\zeta_{\Bbb F}^{k-1}(2\sigma_0)}\sum_{ {\bf n}, {\bf d}, {\bf m}\in \mathcal{O}_{\Bbb F}}\frac{\sigma_{1-\sigma_0}^{k-1}({\bf n}{\bf d})  |\mu({\bf m})|   }
{N^{\sigma_0}({\bf n})   N^{\sigma_0}({\bf m}) N^{2\sigma_0-1}({\bf d})   }  \frac{1}{T\left|\log \frac{X}{N({\bf m}) N({\bf d})  }\right|+1}\nonumber\\
 &&  \leq \frac{\zeta_{\Bbb F}^{k-1}(\sigma_0) }{\zeta_{\Bbb F}^{k-1}(2\sigma_0)}
\sum_{ {\bf n}, {\bf d}, {\bf m}\in \mathcal{O}_{\Bbb F}}\frac{\sigma_{1-\sigma_0}^{k-1}({\bf n} ) \sigma_{1-\sigma_0}^{k-1}({\bf d} ) |\mu({\bf m})|   }
{N^{\sigma_0}({\bf n})   N^{\sigma_0}({\bf m}) N^{2\sigma_0-1}({\bf d})   }  \frac{ \zeta_{\Bbb F}^{k-1}(\sigma_0)}{T\left|\log \frac{X}{N({\bf m}) N({\bf d})  }\right|+1}\nonumber\\
&& \leq \frac{\zeta_{\Bbb F}^{k-1}(\sigma_0) }{\zeta_{\Bbb F}^{k-1}(2\sigma_0)}
\sum_{ {\bf n}, {\bf d}, {\bf m}\in \mathcal{O}_{\Bbb F}}\frac{\sigma_{ 0}^{k-1}({\bf n} ) \sigma_{ 0}^{k-1}({\bf d} ) |\mu({\bf m})|   }
{N^{\sigma_0}({\bf n})   N^{\sigma_0}({\bf m}) N^{\sigma_0}({\bf d})   }  \frac{1}{T\left|\log \frac{X}{N({\bf m}) N({\bf d})  }\right|+1}\nonumber\\
&&=\frac{\zeta_{\Bbb F}^{k-1}(\sigma_0) }{\zeta_{\Bbb F}^{k-1}(2\sigma_0)}
\sum_{ {\bf n} \in \mathcal{O}_{\Bbb F}}\frac{\sigma_{ 0}^{k-1}({\bf n} )   }
{N^{\sigma_0}({\bf n})    }
\sum_{   {\bf d}, {\bf m}\in \mathcal{O}_{\Bbb F}}\frac{ \sigma_{ 0}^{k-1}({\bf d} ) |\mu({\bf m})|   }
{    N^{\sigma_0}({\bf m}) N^{\sigma_0}({\bf d})   }  \frac{1}{T\left|\log \frac{X}{N({\bf m}) N({\bf d})  }\right|+1}\nonumber\\
&&\ll  \zeta_{\Bbb F}^{k-1+2^{k-1}}(\sigma_0)
\sum_{   {\bf d}, {\bf m}\in \mathcal{O}_{\Bbb F}}\frac{ \sigma_{ 0}^{k-1}({\bf d} ) |\mu({\bf m})|   }
{    N^{\sigma_0}({\bf m}) N^{\sigma_0}({\bf d})   }  \frac{1}{T\left|\log \frac{X}{N({\bf m}) N({\bf d})  }\right|+1},\nonumber
\end{eqnarray}
where in the fifth line we used the bound $\sigma_0({\bf d}{\bf n})\leq  \sigma_0({\bf d} ) \sigma_0( {\bf n}) $ and in the final line we used (4.8), which holds for any $k\geq 0.$

Define
$${\bf g}({\bf n}):=\sum_{{\bf n}={\bf m}{\bf d}} \sigma_{ 0}^{k-1}({\bf d} ) |\mu({\bf m})|,\ \ \
G_2(s)=\sum_{   {\bf n} \in \mathcal{O}_{\Bbb F}}\frac{{\bf g}({\bf n})}{ N^s({\bf n}) }\ (\Re s>1). $$
Then we have
\begin{equation}
\sum_{   {\bf d}, {\bf m}\in \mathcal{O}_{\Bbb F}}\frac{ \sigma_{ 0}^{k-1}({\bf d} ) |\mu({\bf m})|   }
{    N^{\sigma_0}({\bf m}) N^{\sigma_0}({\bf d})   }  \frac{1}{T\left|\log \frac{X}{N({\bf m}) N({\bf d})  }\right|+1}=
\sum_{{\bf n}\in \mathcal{O}_{\Bbb F}}\frac{{\bf g}({\bf n})}{N^{\sigma_0}({\bf n})}
\frac{1}{T\left|\log \frac{X}{N({\bf n})    }\right|+1}.
\end{equation}

It is easily seen that ${\bf g}({\bf n})$ is multiplicative. For any prime ideal ${\bf p},$ we have
$${\bf g}({\bf p})=\sigma^{k-1}({\bf p})+|\mu({\bf p})|=2^{k-1}+1.$$
Thus for $\Re(s)>1$ we have
\begin{eqnarray*}
G_2(s)=\sum_{{\bf n}\in \mathcal{O}_{\Bbb F}}\frac{{\bf g}({\bf n})}{N^{s}({\bf n})} =\prod_{{\bf p} \in \mathcal{O}_{\Bbb F}  }
\left(1+\sum_{\alpha=1}^\infty  \frac{{\bf g}({\bf p}^\alpha)}{N^{s}({\bf p}^\alpha)}  \right)
=(\zeta_{\Bbb F}(s))^{2^{k-1}+1}{\bf H}(s),
\end{eqnarray*}
where
 $ {\bf H}(s)$
is absolutely convergent for $\Re(s)>1/2.$ So we have
\begin{equation}
G_2(\sigma_0)\ll (\zeta_{\Bbb F}( \sigma_0))^{2^{k-1}+1}.
\end{equation}

If we write
$$ G_2(s)=\sum_{ n=1}\frac{g_2(n)}{ n^s }, $$
then it is easy to see that $g_2(n)\ll n^\varepsilon.$
So from
 (4.4), (4.10), (4.11),(4.12) and (2.3)   we get
\begin{eqnarray}
 E_{1,k}(X,T)\ll    \frac{X^\varepsilon}{T}.
\end{eqnarray}

\section{\bf Proof of Theorem 1}


Without loss of generality, we suppose that both $X$ and $Y$ are half integers with
$3\leq X<Y$.   Let $  T\geq 3 $ be a parameter to be determined later.
Define
$$
b:=1+\frac{1}{\log X}, \ \  T_1:=T , \ \ T_2:=2 T.
$$

By the definition of $C_{\Bbb F,1}(X,Y)$ and Lemma 2.6 we have
\begin{equation}
C_{\Bbb F,1}(X,Y)=I_{\Bbb F,1}(X,Y,T)+O(XYE_{\Bbb F,1}(X,T)+XY\mathfrak{E}_1(Y,T)),
\end{equation}
where
\begin{eqnarray*}
&&I_{\Bbb F,1}(X,Y,T):=\frac{1}{(2\pi i)^2}\int_{b-iT_1}^{b+iT_1}ds\int_{b-iT_2}^{b+iT_2}
\frac{\zeta_{\Bbb F}(w)\zeta_{\Bbb F}(w+s-1)}{\zeta_{\Bbb F}(s)}\frac{X^sY^w}{sw}dw,\\
&&E_{\Bbb F,1}(X,T):=\sum_{\bf m}\sum_{\bf n}
\frac{|c_{\bf m}({\bf n})|}{N^b({\bf m})N^b({\bf n})}\times \frac{1}{T\left|\frac{X}{N({\bf m})}\right|+1},\\
&&\mathfrak{E}_1(X,T):=\sum_{\bf m}\sum_{\bf n}
\frac{|c_{\bf m}({\bf n})|}{N^b({\bf m})N^b({\bf n})}\times \frac{1}{T\left|\frac{Y}{N({\bf n})}\right|+1}.
\end{eqnarray*}

From (4.9) and (4.13) with $k=1,$ we have
\begin{equation}
E_{\Bbb F,1}(X,T)\ll \frac{X^\varepsilon}{T},\ \ \
\mathfrak{E}_1(X,T)\ll \frac{Y^\varepsilon}{T}.
\end{equation}

We consider the rectangle domain of $w$ formed by the four points
$b\pm iT_2$ and $1/2\pm iT_2.$ Let
$$G(w;x,X,Y):=\frac{\zeta_{\Bbb F}(w)\zeta_{\Bbb F}(w+s-1)}{\zeta_{\Bbb F}(s)}\frac{X^sY^w}{sw}.$$
In this domain,  $G(w;x,X,Y)$
has two simple poles, which are $w=1$ and $w=2-s,$ respectively. It is easy to
see that
\begin{eqnarray*}
&&Res_{w=1}G(w;s,X,Y)=\rho_{\Bbb F}Y\frac{X^s}{s},
\\&&
Res_{w=2-s}G(w;s,X,Y)= \rho_{\Bbb F}\frac{\zeta_{\Bbb F}(2-s)}{\zeta_{\Bbb F}(s)}
\frac{X^sY^{2-s}}{s(2-s)}.
\end{eqnarray*}
By the residue theorem we get
\begin{eqnarray}
I_{\Bbb F,1}(X,Y,T)&&=\mathfrak{J}_{1}(X,Y,T)+\mathfrak{J}_{2}(X,Y,T)\\
&&+H_1(X,Y,T)+H_2(X,Y,T)-H_3(X,Y,T),\nonumber
\end{eqnarray}
where
\begin{eqnarray}
&&\mathfrak{J}_{1}(X,Y,T):=\rho_{\Bbb F}Y\frac{1}{ 2\pi i }\int_{b-iT_1}^{b+iT_1}\frac{X^s}{s}ds,\\
&&\mathfrak{J}_{2}(X,Y,T):=\rho_{\Bbb F}\frac{1}{ 2\pi i }\int_{b-iT_1}^{b+iT_1}
 \frac{\zeta_{\Bbb F}(2-s)}{\zeta_{\Bbb F}(s)}
\frac{X^sY^{2-s}}{s(2-s)}ds,\nonumber\\
&&H_1(X,Y,T):=\frac{1}{(2\pi i)^2}\int_{b-iT_1}^{b+iT_1}ds\int_{1/2+iT_2}^{b+iT_2}
\frac{\zeta_{\Bbb F}(w)\zeta_{\Bbb F}(w+s-1)}{\zeta_{\Bbb F}(s)}\frac{X^sY^w}{sw}dw,\nonumber\\
&&H_2(X,Y,T):=\frac{1}{(2\pi i)^2}\int_{b-iT_1}^{b+iT_1}ds\int_{1/2-iT_2}^{1/2+iT_2}
\frac{\zeta_{\Bbb F}(w)\zeta_{\Bbb F}(w+s-1)}{\zeta_{\Bbb F}(s)}\frac{X^sY^w}{sw}dw,\nonumber\\
&&H_3(X,Y,T):=\frac{1}{(2\pi i)^2}\int_{b-iT_1}^{b+iT_1}ds\int_{1/2-iT_2}^{b-iT_2}
\frac{\zeta_{\Bbb F}(w)\zeta_{\Bbb F}(w+s-1)}{\zeta_{\Bbb F}(s)}\frac{X^sY^w}{sw}dw.\nonumber
\end{eqnarray}

We consider $H_1(X,Y,T)$ first. Suppose
$$s=b+it,\ |t|\leq T, \ \ w=u+2iT,\ 1/2\leq u\leq b.$$
From Lemma 2.2 we have
\begin{eqnarray*}
G(w;x,X,Y)\ll\left\{\begin{array}{ll}
\frac{X^bY^uT^{1-2u}}{|t|+1}\log^3 T,&\mbox{$1/2\leq u\leq 1,$} \\
\frac{X^bY^u}{(|t|+1)T}\log^3 T ,&\mbox{$1\leq u\leq b,$}
\end{array}\right.
\end{eqnarray*}
which implies that
\begin{eqnarray}
H_1(X,Y,T)&&\ll \int_{b-iT_1}^{b+iT_1}\frac{X^b}{|t|+1}dt
\left(\int_{1/2}^1 Y^uT^{1-2u} du +\int_{1}^b \frac{Y^u}{T}  du \right)\log^3 T\\
&&\ll XY^{1/2}\log^4 T+XY^bT^{-1}\log^4 T.\nonumber
\end{eqnarray}

Similarly we have
\begin{eqnarray}
H_3(X,Y,T) \ll   XY^{1/2}\log^4 T+XY^bT^{-1}\log^4 T.
\end{eqnarray}

Now we consider $H_2(X,Y,T).$ Suppose
$$s=b+it,\ |t|\leq T,\  w=1/2+iv, \ |v|\leq 2T.$$
We have
\begin{eqnarray}
H_2(X,Y,T) \ll   XY^{1/2}\log T \times \mathfrak{H}(X,Y,T),
\end{eqnarray}
where
\begin{eqnarray*}
\ \ \ \ \ \ \mathfrak{H}(X,Y,T):=
 \int_{b-iT_1}^{b+iT_1}dt\int_{1/2-iT_2}^{1/2+iT_2}
\frac{|\zeta_{\Bbb F}(\frac 12+iv)\zeta_{\Bbb F}(\frac 12+\frac{1}{\log X}+i(t+v))|}
{(|t|+1)(|v|+1) }dv.
\end{eqnarray*}
Write
\begin{eqnarray}
  \mathfrak{H}(X,Y,T)=\mathfrak{H}_1(X,Y,T)+\mathfrak{H}_2(X,Y,T),
\end{eqnarray}
where
\begin{eqnarray*}
&&  \mathfrak{H}_1(X,Y,T):=
 \int_{|t|\leq |v|}
\frac{|\zeta_{\Bbb F}(\frac 12+iv)\zeta_{\Bbb F}(\frac 12+\frac{1}{\log X}+i(t+v))|}
{(|t|+1)(|v|+1) }dvdt,\\
&&  \mathfrak{H}_2(X,Y,T):=
 \int_{|v|\leq |t|}
\frac{|\zeta_{\Bbb F}(\frac 12+iv)\zeta_{\Bbb F}(\frac 12+\frac{1}{\log X}+i(t+v))|}
{(|t|+1)(|v|+1) }dvdt.
\end{eqnarray*}
From Lemma  2.3 and partial integration we have that
\begin{equation}
\int_{-U}^U\frac{|\zeta_{\Bbb F}(u+iv)|^2}{|v|+1}dv\ll (\log U)^5
\ \ (1/2\leq u\leq 1).
\end{equation}
and
\begin{equation}
\int_{-U}^U\frac{|\zeta_{\Bbb F}(u+iv)|}{|v|+1}dv\ll (\log U)^3 \ \ (1/2\leq u\leq 1).
\end{equation}
If $|t|\leq |v|,$ then $|v+t|\leq |v|+|t|\leq 2|v|$, which combining with
(5.9) and Cauchy's inequality implies that
\begin{eqnarray}
&&\ \ \ \ \ \ \ \ \ \ \ \mathfrak{H}_1(X,Y,T)\\&&\ll
 \int_{|t|\leq T}\frac{1}{|t|+1}dt\int_{|t|\leq |v|}
\frac{|\zeta_{\Bbb F}(\frac 12+iv)|}{(|v|+1)^{1/2} }
\frac{|\zeta_{\Bbb F}(\frac 12+\frac{1}{\log X}+i(t+v))|}
{ (|v+t|+1)^{1/2} }dv\nonumber\\
&&\ll \int_{|t|\leq T}\frac{dt}{|t|+1}
\left(\int_{|t|\leq |v|}\frac{|\zeta_{\Bbb F}(\frac 12+iv)|^2}{|v|+1  }dv\right)^{\frac 12}
\left(\int_{|t|\leq |v|}\frac{|\zeta_{\Bbb F}(\frac 12+\frac{1}{\log X}+i(v+t))|^2}{|v+t|+1  }dv\right)^{\frac 12}\nonumber\\
&&\ll (\log T)^6.\nonumber
\end{eqnarray}

If $|v|\leq |t|,$ then $|v+t|\leq |v|+|t|\leq 2|t|$, which combining with (5.10) gives
\begin{eqnarray}
\ \ \ \ \mathfrak{H}_2(X,Y,T)&& \ll
 \int_{|v|\leq T} \frac{|\zeta_{\Bbb F}(\frac 12+iv)|dv}{|v|+1 }
 \int_{|v|\leq |t|}
\frac{|\zeta_{\Bbb F}(\frac 12+\frac{1}{\log X}+i(t+v))|dt}
{  |v+t|+1  }\\
&&\ll (\log T)^6.\nonumber
\end{eqnarray}

From (5.7), (5.8), (5.11) and (5.12) we get
\begin{eqnarray}
H_2(X,Y,T) \ll   XY^{1/2}(\log T)^7.
\end{eqnarray}

From (2.15) we get
\begin{eqnarray}
&&\mathfrak{J}_{1}(X,Y,T)=\rho_{\Bbb F}Y +O\left(\frac{XY}{T\log X}\right).
\end{eqnarray}

Finally we consider $\mathfrak{J}_2(X,Y,T).$
Let
\begin{eqnarray*}
&&IP_1=\{s=\sigma-iT_1: b\leq \sigma\leq 2\},\ \
IP_2=\{s=2+it: -T\leq t\leq -\frac{1}{\log X}\},\\
&&IP_3=\{s=\frac{e^{i\theta}}{\log X}: -\frac{\pi}{2}\leq \theta\leq \frac{\pi}{2}\},\\
&&IP_4=\{s=2+it:  \frac{1}{\log X}\leq t\leq T\},\ \
IP_5=\{s=\sigma+iT_1: b\leq \sigma\leq 2\}.
\end{eqnarray*}

By the residue theorem we have
\begin{equation}
\mathfrak{J}_2(X,Y,T)=\frac{\zeta_{\Bbb F}(0)}{2\zeta_{\Bbb F}(2)}X^2+
 \sum_{j=1}^4\mathcal{J}_{2j}(X,Y,T)-\mathcal{J}_{25}(X,Y,T),
\end{equation}
where
\begin{eqnarray*}
&&\mathfrak{J}_{2j}(X,Y,T):=\rho_{\Bbb F}\frac{1}{ 2\pi i }
\int_{IP_j } \frac{\zeta_{\Bbb F}(2-s)}{\zeta_{\Bbb F}(s)}
\frac{X^sY^{2-s}}{s(2-s)}ds\ \ (j=1, 2, 3, 4, 5).
\end{eqnarray*}

By Lemma 2.2 we have
 \begin{eqnarray}
\ \ \mathfrak{J}_{21}(X,Y,T) \ll \int_b^{2 }
T^{\sigma-3}X^\sigma Y^{2 -\sigma}  \log^2 Td\sigma
 \ll \frac{XY}{T^2} \log^2 T+\frac{X^2}{T}\log^2 T.
\end{eqnarray}
and
 \begin{eqnarray}
\ \ \mathfrak{J}_{25}(X,Y,T) \ll \int_b^{2 }
T^{\sigma-3}X^\sigma Y^{2 -\sigma}  \log^2 Td\sigma
 \ll \frac{XY}{T^2} \log^2 T+\frac{X^2}{T}\log^2 T.
\end{eqnarray}

 By Lemma 2.2 again we have
\begin{eqnarray}
\ \ \ \  \mathfrak{J}_{24}(X,Y,T) &&\ll   X^{2}
 \int_{\frac{1}{\log X}}^{T_1}\frac{\log(t+1)}{t(t+1)}dt\\
 &&=X^{2}
 \left(\int_{\frac{1}{\log X}}^{1}\frac{\log(t+1)}{t(t+1)}dt+
 \int_{1}^{T_1}\frac{\log(t+1)}{t(t+1)}dt
 \right)
 \ll X^2.\nonumber
\end{eqnarray}
by noting that $\log(1+t)\ll t\ (0<t<1).$
Similarly
\begin{eqnarray}
 &&\ \  \ \ \  \mathfrak{J}_{22}(X,Y,T) \ll   X^{2}
 \int_{\frac{1}{\log X}}^{T_1}\frac{\log(t+1)}{t(t+1)}dt
 \ll X^2.
\end{eqnarray}

For $\mathfrak{J}_{23}(X,Y,T)$ we have
 \begin{eqnarray}
 &&\ \  \ \ \  \mathfrak{J}_{23}(X,Y,T) \ll   X^{2}.
\end{eqnarray}

From (5.15)-(5.20) we get
\begin{equation}
\mathfrak{J}_2(X,Y,T)\ll  \frac{XY}{T^2} \log^2 T
+\frac{X^2}{T}\log^2 T+X^{2}.
\end{equation}

From (5.1)-(5.6),   (5.13), (5.14) and (5.21) we have
\begin{eqnarray}
C_{\Bbb F,1}(X, Y)&& =\rho_{\Bbb F}Y+O\left(\frac{XY^{1+\varepsilon}}{T} +XY^{1/2}\log^7 T\right)\\
&&\ \ \ \ \ \ \ \ \ \ +O\left(X^{2 }
+\frac{X^2}{T}\log^4 T\right)\nonumber\\
&&=\rho_{\Bbb F}Y+O\left(XY^{1/2} \log^7 Y +X^{2 } \right)\nonumber
\end{eqnarray}
by choosing $T=XY.$ This completes the proof of Theorem 1.

\section{\bf Proof of Theorem 2}

Without loss of generality, we suppose that both $X$ and $Y$ are half integers and
$X<Y$.   Let $ T \geq 3$ be a parameter to be determined later.
Define
\begin{eqnarray*}
&&b_1:=1+\frac{1}{\log X},\ b_2:=1+\frac{2}{\log X},\ b_3:=1+\frac{3}{\log X} \\
&&  T_1=T, \ \ T_2=2T, \ \ T_3:=4T.
\end{eqnarray*}

By the definition of $C_{\Bbb F,2}(X,Y)$ and Lemma 2.6 we have
\begin{equation}
C_{\Bbb F,2}(X,Y)=I_{\Bbb F,2}(X,Y,T)+O(X^2Y^{b}E_{\Bbb F,1}(X,T)+X^2Y^{b}\mathfrak{E}_2(Y,T)),
\end{equation}
where
\begin{eqnarray*}
&&I_{\Bbb F,2}(X,Y,T):=\frac{1}{(2\pi i)^3}\int_{b_1-iT_1}^{b_1+iT_1}ds_1\int_{b_2-iT_2}^{b_2+iT_2}ds_2
\int_{b_3-iT_3}^{b_3+iT_3}
 \mathcal{G}(w;s_1,s_2)dw,\\
&&E_{\Bbb F,1}(X,T):=\sum_{\bf m_1}\sum_{\bf m_2}\sum_{\bf n}
\frac{|c_{\bf m_1}({\bf n})c_{\bf m_2}({\bf n})|}{N^{b_1}({\bf m_1})N^{b_1}({\bf m_2})N^{b_1}({\bf n})}\times \frac{1}{T\left|\frac{X}{N({\bf m})}\right|+1},\\
&&\mathfrak{E}_2(X,T):=\sum_{\bf m_1}\sum_{\bf m_2}\sum_{\bf n}
\frac{|c_{\bf m_1}({\bf n})c_{\bf m_2}({\bf n})|}{N^{b_1}({\bf m_1})N^{b_1}({\bf m_2})N^{b_1}({\bf n})}\times \frac{1}{T\left|\frac{Y}{N({\bf n})}\right|+1}
\end{eqnarray*}
and
$$\mathcal{G}(w;s_1,s_2):=\frac{ \zeta_{\Bbb F}(w)   \zeta_{\Bbb F}(w+s_1-1)  \zeta_{\Bbb F}(w+s_2-1)  \zeta_{\Bbb F}(w+s_1+s_2-2)    }
{\zeta_{\Bbb F}(s_1)\zeta_{\Bbb F}(s_2)  \zeta_{\Bbb F}(2w+s_1+s_2-2)}
\frac{X^{s_1+s_2}Y^{w}}{s_1s_2w}.
$$

From (4.9) and (4.13) with $k=2,$ we have
\begin{equation}
E_{\Bbb F,1}(X,T)\ll \frac{X^\varepsilon}{T},\ \ \
\mathfrak{E}_2(X,T)\ll \frac{Y^\varepsilon}{T}.
\end{equation}

We consider the rectangle domain of $w$ formed by the four points
$b_3\pm iT_3$ and $2/3\pm iT_3.$
In this domain,  $\mathcal{G}(w;s_1,s_2) $
has four simple poles, which are $w_1=1$,   $w_2=2-s_1,$  $w_3=2-s_2$ and $ w_4=3-s_1-s_2$  respectively.
By the residue theorem we get
\begin{eqnarray}
I_{\Bbb F,2}(X,Y,T)&&=\mathfrak{L}_{1}(X,Y,T)+\mathfrak{L}_{2}(X,Y,T)
+\mathfrak{L}_{3}(X,Y,T)+\mathfrak{L}_{4}(X,Y,T)
\\
&&+K_1(X,Y,T)+K_2(X,Y,T)-K_3(X,Y,T),\nonumber
\end{eqnarray}
where
\begin{eqnarray*}
&&\mathfrak{L}_{j}(X,Y,T):= \frac{1}{ (2\pi i)^2 }\int_{b_1-iT_1}^{b_1+iT_1}ds_1\int_{b_2-iT_2}^{b_2+iT_2}Res_{w=w_j}\mathcal{G}(w;s_1,s_2) ds_2\ \ (j=1,2,3,4),\\
&&K_1(X,Y,T):=\frac{1}{(2\pi i)^3}\int_{b_1-iT_1}^{b_1+iT_1}ds_1\int_{b_2-iT_2}^{b_2+iT_2}ds_2
\int_{2/3+iT_3}^{b_3+iT_3}\mathcal{G}(w;s_1,s_2)
 dw,\nonumber\\
&&K_2(X,Y,T):=\frac{1}{(2\pi i)^3}\int_{b_1-iT_1}^{b_1+iT_1}ds_1\int_{b_2-iT_2}^{b_2+iT_2}ds_2
\int_{2/3-iT_3}^{2/3+iT_3}\mathcal{G}(w;s_1,s_2)
 dw,\nonumber\\
&&K_3(X,Y,T):=\frac{1}{(2\pi i)^3}\int_{b_1-iT_1}^{b_1+iT_1}ds_1\int_{b_2-iT_2}^{b_2+iT_2}ds_2
\int_{2/3-iT_3}^{b_3-iT_3}\mathcal{G}(w;s_1,s_2) dw.\nonumber
\end{eqnarray*}

\subsection{\bf Upper bound of $K_j(X,Y,T)$(j=1,2,3)}\
We first consider $K_1(X,Y,T).$
Suppose $$s_1=b+it_1, |t_1|\leq T_1,  s_2=b+it_2, |t_2|\leq T_2, w=u+iT_3, 2/3\leq u\leq b.$$
By Lemma 2.5 we have
$$
\mathcal{G}(w;s_1,s_2)\ll
\frac{X^{2b}}{(|t_1|+1)(|t_2|+1)}
T^{\frac 53-\frac{8u}{3} }Y^{u}\log^6 T,\ \ 2/3\leq u\leq b,
$$
which implies that
\begin{equation}
K_1(X,Y,T)\ll  X^2 Y^{\frac 23}T^{-\frac 19}\log^8 T+X^2Y^bT^{-1}\log^8 T.
\end{equation}

Similarly we have
\begin{equation}
K_3(X,Y,T)\ll  X^2 Y^{\frac 23}T^{-\frac 19}\log^8 T+X^2Y^bT^{-1}\log^8 T.
\end{equation}

We now consider $K_2(X,Y,T).$ By (2.4) of  Lemma 2.2   we write
\begin{equation}
K_2(X,Y,T)\ll \mathfrak{J}X^2 Y^{\frac 23}\log^2 T ,
\end{equation}
where
\begin{equation}
\mathfrak{J} :=
 \int_{b_1-iT_1}^{b_1+iT_1}dt_1 \int_{b_2-iT_2}^{b_2+iT_2}dt_2
 \int_{2/3-iT_3}^{2/3+iT_3}
\frac{|g(t_1,t_2,v)|}
{(|t_1|+1)(|t_2|+1)(|v|+1) }dv
\end{equation}
 with
 \begin{eqnarray*}
g(t_1,t_2,v):&&=\zeta_{\Bbb F}\left(\frac 23+iv\right)
\zeta_{\Bbb F}\left(\frac 23+\frac{1}{\log X}+i(v+t_1)\right)\\
&&\ \ \ \ \   \zeta_{\Bbb F}\left(\frac 23+\frac{1}{\log X}+i(v+t_2)\right)
 \zeta_{\Bbb F}\left(\frac 23+\frac{2}{\log X}+i(v+t_1+t_2)\right)
\end{eqnarray*}

With the help of  Lemma 2.4  we can show that
\begin{eqnarray}
\mathfrak{J} \ll \log^3 T.
\end{eqnarray}
The proof of (6.8) is similar to  the arguments of $H_2(T)$ in T\'{o}th and Zhai
\cite{TZ}. So we omit its details.

From (6.6)-(6.8)  we get
\begin{equation}
K_2(X,Y,T)\ll  X^2 Y^{\frac 23} \log^5 T.
\end{equation}

\subsection{\bf Evaluaton of $\mathfrak{L}_{1}(X,Y,T) $ }\
It is easy to see that
$$Res_{w=1}\mathcal{G}(w;s_1,s_2) =\rho_{\Bbb F}\frac{\zeta_{\Bbb F}(s_1+s_2-1)    }
{  \zeta_{\Bbb F}(s_1+s_2)}
\frac{X^{s_1+s_2}Y }{s_1s_2}. $$

So
$$\mathfrak{L}_{1}(X,Y,T)=\frac{1}{ (2\pi i)^2 }\int_{b_1-iT_1}^{b_1+iT_1}ds_1\int_{b_2-iT_2}^{b_2+iT_2}
\rho_{\Bbb F}\frac{\zeta_{\Bbb F}(s_1+s_2-1)    }
{  \zeta_{\Bbb F}(s_1+s_2)}
\frac{X^{s_1+s_2}Y }{s_1s_2} ds_2.
$$

We consider the rectangle domain of $s_2$ formed by the four points $1/2\pm iT_2 $ and $b\pm iT_2.$ In this domain, the integral function $\rho_{\Bbb F}\frac{\zeta_{\Bbb F}(s_1+s_2-1)    }
{  \zeta_{\Bbb F}(s_1+s_2)}
\frac{X^{s_1+s_2}Y }{s_1s_2} $ has a simple pole $s_2=2-s_1$,
with residue $\frac{\rho_{\Bbb F}^2}{\zeta_{\Bbb F}(2)}\frac{X^2Y}{s_1(2-s_1)}.$
By the residue theorem, we have
\begin{eqnarray}
\mathfrak{L}_{1}(X,Y,T)&&=\frac{1}{2\pi i}\int_{b_1-iT_1}^{b_1+iT_1}\frac{\rho_{\Bbb F}^2}{\zeta_{\Bbb F}(2)}\frac{X^2Y}{s_1(2-s_1)}ds_1\\
&&\ \ \ +\mathfrak{L}_{11}(X,Y,T)+\mathfrak{L}_{12}(X,Y,T)-\mathfrak{L}_{13}(X,Y,T)
\nonumber
\end{eqnarray}
where
\begin{eqnarray*}
&&\mathfrak{L}_{11}(X,Y,T):=\frac{1}{ (2\pi i)^2 }\int_{b_1-iT_1}^{b_1+iT_1}ds_1\int_{1/2+iT_2}^{b_2+iT_2}\rho_{\Bbb F}\frac{\zeta_{\Bbb F}(s_1+s_2-1)    }
{  \zeta_{\Bbb F}(s_1+s_2)}
\frac{X^{s_1+s_2}Y }{s_1s_2} ds_2,\\
&&\mathfrak{L}_{12}(X,Y,T):=\frac{1}{ (2\pi i)^2 }\int_{b_1-iT_1}^{b_1+iT_1}ds_1\int_{1/2-iT_2}^{1/2+iT_2}\rho_{\Bbb F}\frac{\zeta_{\Bbb F}(s_1+s_2-1)    }
{  \zeta_{\Bbb F}(s_1+s_2)}
\frac{X^{s_1+s_2}Y }{s_1s_2} ds_2,\\
&&\mathfrak{L}_{13}(X,Y,T):=\frac{1}{ (2\pi i)^2 }\int_{b_1-iT_1}^{b_1+iT_1}ds_1\int_{1/2-iT_2}^{b_2-iT_2}\rho_{\Bbb F}\frac{\zeta_{\Bbb F}(s_1+s_2-1)    }
{  \zeta_{\Bbb F}(s_1+s_2)}
\frac{X^{s_1+s_2}Y }{s_1s_2} ds_2.
\end{eqnarray*}

By Lemma 2.2 it is easy to see that
\begin{eqnarray}
&&\mathfrak{L}_{11}(X,Y,T)\ll   \frac{X^2Y}{T}\log T+\frac{X^{\frac 32}Y}{T^{\frac 12}}\log T
\end{eqnarray}
and
\begin{eqnarray}
&&\mathfrak{L}_{12}(X,Y,T)\ll  \frac{X^2Y}{T}\log T+\frac{X^{\frac 32}Y}{T^{\frac 12}}\log T.
\end{eqnarray}

For $\mathfrak{L}_{12}(X,Y,T)$,   we have
\begin{eqnarray}
\mathfrak{L}_{12}(X,Y,T)&&\ll YX^{\frac 32}\int_{-T_1}^{T_1}\frac{dt_1}{|t_1|+1}
\int_{-T_2}^{T_2}\frac{|\zeta(\frac 12+\frac{1}{\log X}+i(t_1+t_2))|}{|t_2|+1}dt_2\\
&&\ll YX^{\frac 32}(\log T)^4,\nonumber
\end{eqnarray}
where we used the estimate
\begin{eqnarray*}
 \int_{-T_1}^{T_1}\frac{dt_1}{|t_1|+1}
\int_{-T_2}^{T_2}\frac{|\zeta(\frac 12+\frac{1}{\log X}+i(t_1+t_2))|}{|t_2|+1}dt_2
 \ll  (\log T)^4,
\end{eqnarray*}
whose proof is similar to (5.11) and (5.12).

We now consider the first integral in (6.10). We have
\begin{eqnarray*}
 \frac{1}{2\pi i}\int_{b-iT_1}^{b+iT_1}\frac{\rho_{\Bbb F}^2}{\zeta_{\Bbb F}(2)}\frac{X^2Y}{s_1(2-s_1)}ds_1=
 \frac{1}{2\pi i}\int_{(b)} \frac{\rho_{\Bbb F}^2}{\zeta_{\Bbb F}(2)}\frac{X^2Y}{s_1(2-s_1)}ds_1+O\left(\frac{X^2Y\log T}{T}\right),
\end{eqnarray*}
where
$\int_{(b)}$ means that $\int_{b-i\infty}^{b+i\infty}.$ Moving the integral line from
$b$ to $\Re(s_1)=1,$ we get
\begin{eqnarray}
 \frac{1}{2\pi i}\int_{b-iT_1}^{b+iT_1}\frac{\rho_{\Bbb F}^2}{\zeta_{\Bbb F}(2)}\frac{X^2Y}{s_1(2-s_1)}ds_1=c_{\Bbb F}X^2Y
 +O\left(\frac{X^2Y\log T}{T}\right),
\end{eqnarray}
with
\begin{eqnarray*}
c_{\Bbb F}=\frac{\rho_{\Bbb F}^2}{\zeta_{\Bbb F}(2)}\frac{1}{2\pi i}\int_{(1)}\frac{1}{s_1(2-s_1)}ds_1=\frac{\rho_{\Bbb F}^2}{2\zeta_{\Bbb F}(2)}.
\end{eqnarray*}

From (6.10)-(6.14)   we get
\begin{eqnarray}
\mathfrak{L}_{1}(X,Y,T)=\frac{\rho_{\Bbb F}^2}{2\zeta_{\Bbb F}(2)} X^2Y +O\left(\frac{X^2Y\log T}{T}+YX^{\frac 32}(\log T)^4\right).
\end{eqnarray}

\subsection{\bf Upper bound of $\mathfrak{L}_{2}(X,Y,T) $ }\
It is easy to see that
$$Res_{w=2-s_1}\mathcal{G}(w;s_1,s_2) = \rho_{\Bbb F}
\frac{\zeta_{\Bbb F}(2-s_1)  \zeta_{\Bbb F}(1-s_1+s_2)    }
{\zeta_{\Bbb F}(s_1)   \zeta_{\Bbb F}(2-s_1+s_2)}
\frac{X^{s_1+s_2}Y^{2-s_1}}{s_1s_2(2-s_1)}$$
So we have
\begin{eqnarray*}
\mathfrak{L}_{2}(X,Y,T)=\frac{1}{ (2\pi i)^2 }\int_{b_1-iT_1}^{b_1+iT_1}ds_1\int_{b_2-iT_2}^{b_2+iT_2} \rho_{\Bbb F}
\frac{\zeta_{\Bbb F}(2-s_1)  \zeta_{\Bbb F}(1-s_1+s_2)    }
{\zeta_{\Bbb F}(s_1)   \zeta_{\Bbb F}(2-s_1+s_2)}
\frac{X^{s_1+s_2}Y^{2-s_1}}{s_1s_2(2-s_1)} ds_2.
\end{eqnarray*}

We consider the rectangle domain of $s_2$  formed by the four points
$1/2\pm iT_2$ and $b_2\pm iT_2.$ In this domain, the integral function
in the above integral is $s_2=s_1$ with residue
\begin{eqnarray*}
\frac{\rho_{\Bbb F}^2}{\zeta_{\Bbb F}(2)}\cdot\frac{\zeta_{\Bbb F}(2-s_1) }
{\zeta_{\Bbb F}(s_1)}\cdot
\frac{X^{2s_1 }Y^{2-s_1}}{s_1^2(2-s_1)}.
\end{eqnarray*}

By the residue theorem we get
\begin{eqnarray}
\ \ \ \ \ \  \ \ \ \mathfrak{L}_{2}(X,Y,T)=\mathfrak{L}_{20}(X,Y,T)+\mathfrak{L}_{21}(X,Y,T)
+\mathfrak{L}_{22}(X,Y,T)-\mathfrak{L}_{23}(X,Y,T)
\end{eqnarray}
where
\begin{eqnarray*}
&&\mathfrak{L}_{20}(X,Y,T):=\frac{\rho_{\Bbb F}^2}{\zeta_{\Bbb F}(2)}\cdot\frac{1}{ 2\pi i  }\int_{b_1-iT_1}^{b_1+iT_1}
\frac{\zeta_{\Bbb F}(2-s_1) }
{\zeta_{\Bbb F}(s_1)}\cdot
\frac{X^{2s_1 }Y^{2-s_1}}{s_1^2(2-s_1)}ds_1    \\
&&\mathfrak{L}_{21}(X,Y,T):=\frac{1}{ (2\pi i)^2}\int_{b_1-iT_1}^{b_1+iT_1}ds_1\int_{1/2+iT_2}^{b_2+iT_2} \rho_{\Bbb F}
\frac{\zeta_{\Bbb F}(2-s_1)  \zeta_{\Bbb F}(1-s_1+s_2)    }
{\zeta_{\Bbb F}(s_1)   \zeta_{\Bbb F}(2-s_1+s_2)}
\frac{X^{s_1+s_2}Y^{2-s_1}}{s_1s_2(2-s_1)} ds_2,\\
&&\mathfrak{L}_{22}(X,Y,T):=\frac{1}{ (2\pi i)^2 }\int_{b_1-iT_1}^{b_1+iT_1}ds_1\int_{1/2-iT_2}^{1/2+iT_2} \rho_{\Bbb F}
\frac{\zeta_{\Bbb F}(2-s_1)  \zeta_{\Bbb F}(1-s_1+s_2)    }
{\zeta_{\Bbb F}(s_1)   \zeta_{\Bbb F}(2-s_1+s_2)}
\frac{X^{s_1+s_2}Y^{2-s_1}}{s_1s_2(2-s_1)} ds_2,\\
&&\mathfrak{L}_{23}(X,Y,T):=\frac{1}{ (2\pi i)^2 }\int_{b_1-iT_1}^{b_1+iT_1}ds_1\int_{1/2-iT_2}^{b_2-iT_2} \rho_{\Bbb F}
\frac{\zeta_{\Bbb F}(2-s_1)  \zeta_{\Bbb F}(1-s_1+s_2)    }
{\zeta_{\Bbb F}(s_1)   \zeta_{\Bbb F}(2-s_1+s_2)}
\frac{X^{s_1+s_2}Y^{2-s_1}}{s_1s_2(2-s_1)} ds_2.
\end{eqnarray*}

By Lemma 2.2 we have
\begin{eqnarray}
\ \ \ \ \ \ \ \ \  \mathfrak{L}_{21}(X,Y,T)&&\ll Y\log T\int_{-T_1}^{T_1}
\frac{|\zeta_{\Bbb F}(2-b_1-it_1)|    }
{ (|t_1|+1)^2  } dt_1
\int_{1/2 }^{b_2 }  T^{-\sigma_2} X^{b_1+\sigma_2}  d\sigma_2 \\
&&\ll \frac{X^2Y\log T}{T}+\frac{YX^{ 3/2}\log T}{T^{  1/2}}\nonumber
\end{eqnarray}
and
\begin{eqnarray}
\mathfrak{L}_{23}(X,Y,T)\ll  \frac{X^2Y\log T}{T}+\frac{YX^{ 3/2}\log T}{T^{  1/2}}.
\end{eqnarray}

By (2.4) of Lemma 2.2 we have
\begin{eqnarray}
\ \ \ \ \ \mathfrak{L}_{22}(X,Y,T)&&\ll  X^{\frac 32}Y\log^2  T\int_{-T_1}^{ T_1}dt_1\int_{-T_2}^{T_2}
\frac{  | \zeta_{\Bbb F}(\frac 12-\frac{1}{\log X}+i(t_2-t_1))| }
 {(|t_1|+1)^2(|t_2|+1)} dt_2\\
 &&\ll X^{\frac 32}Y\log^5 T,\nonumber
\end{eqnarray}
where we used the bound
$$\int_{-T_1}^{ T_1}dt_1\int_{-T_2}^{T_2}
\frac{  |\zeta_{\Bbb F}(\frac 12-\frac{1}{\log X}+i(t_2-t_1)) |    }
 {(|t_1|+1)^2(|t_2|+1)} dt_2\ll \log^3 T,$$
whose proof is similar to (5.11) and (5.12).

Finally we consider  $\mathfrak{L}_{20}(X,Y,T)$. We write
\begin{eqnarray*}
\mathfrak{L}_{20}(X,Y,T)= \frac{\rho_{\Bbb F}^2}{\zeta_{\Bbb F}(2)}\cdot\frac{1}{ 2\pi i  }\int_{(b_1)}
\frac{\zeta_{\Bbb F}(2-s_1) }
{\zeta_{\Bbb F}(s_1)}\cdot
\frac{X^{2s_1 }Y^{2-s_1}}{s_1^2(2-s_1)}ds_1+O\left(\frac{X^2Y\log^2 T}{T^2}\right).
\end{eqnarray*}
Moving the integral line to $\Re s_1=12/5$, we encounter a simple pole $s_1=2.$  We have
\begin{eqnarray*}
&& \ \ \ \ \ \ \frac{\rho_{\Bbb F}^2}{\zeta_{\Bbb F}(2)}\cdot\frac{1}{ 2\pi i  }\int_{(b_1)}
\frac{\zeta_{\Bbb F}(2-s_1) }
{\zeta_{\Bbb F}(s_1)}\cdot
\frac{X^{2s_1 }Y^{2-s_1}}{s_1^2(2-s_1)}ds_1\\&&=\frac{\zeta_{\Bbb F}(0)\rho_{\Bbb F}^2 }
{4\zeta_{\Bbb F}^2(2)}X^4
+\frac{\rho_{\Bbb F}^2}{\zeta_{\Bbb F}(2)}\cdot\frac{1}{ 2\pi i  }\int_{(12/5)}
\frac{\zeta_{\Bbb F}(2-s_1) }
{\zeta_{\Bbb F}(s_1)}\cdot
\frac{X^{2s_1 }Y^{2-s_1}}{s_1^2(2-s_1)}ds_1.
\end{eqnarray*}

By Lemma 2.1 we see that if $s_1=12/5+it,$ then
$$\frac{\zeta_{\Bbb F}(2-s_1)}{\zeta_{\Bbb F}(s_1)s_1^2(2-s_1)}\ll \frac{1}{(|t|+1)^{6/5}},$$
which implies that
\begin{eqnarray*}
 \int_{(12/5)}
\frac{\zeta_{\Bbb F}(2-s_1) }
{\zeta_{\Bbb F}(s_1)}\cdot
\frac{X^{2s_1 }Y^{2-s_1}}{s_1^2(2-s_1)}ds_1\ll X^{\frac{24}{5}}Y^{-\frac 25}.
\end{eqnarray*}

From the above estimates we have
\begin{eqnarray}
\ \ \ \ \ \mathfrak{L}_{20}(X,Y,T)
&&=\frac{\zeta_{\Bbb F}(0)\rho_{\Bbb F}^2 }
{4\zeta_{\Bbb F}^2(2)}X^4+O\left(\frac{X^2Y\log^2 T}{T^2}+ X^{\frac{24}{5}}Y^{-\frac 25}\right).
\end{eqnarray}

From (6.16) to (6.20) we get
\begin{eqnarray}
\ \ \ \ \ \ \ \ \mathfrak{L}_{2}(X,Y,T)=\frac{\zeta_{\Bbb F}(0)\rho_{\Bbb F}^2 }
{4\zeta_{\Bbb F}^2(2)}X^4 +O\left(\frac{X^2Y\log T}{T}+X^{\frac 32}Y\log^5 T+
X^{\frac{24}{5}}Y^{-\frac 25} \right).
\end{eqnarray}

\subsection{\bf Upper bound of $\mathfrak{L}_{3}(X,Y,T) $ }\
It is easy to see that
\begin{eqnarray*}
&&Res_{w=2-s_2}\mathcal{G}(w;s_1,s_2) =\rho_{\Bbb F}
\frac{\zeta_{\Bbb F}(2-s_2)  \zeta_{\Bbb F}(1+s_1-s_2)    }
{\zeta_{\Bbb F}(s_2)   \zeta_{\Bbb F}(2-s_1+s_2)}
\frac{X^{s_1+s_2}Y^{2-s_2}}{s_1s_2(2-s_2)}.  \nonumber
\end{eqnarray*}
So we have
\begin{eqnarray*}
\mathfrak{L}_{2}(X,Y,T)=\frac{1}{ (2\pi i)^2 }\int_{b_1-iT_1}^{b_1+iT_1}ds_1\int_{b_2-iT_2}^{b_2+iT_2}
\rho_{\Bbb F}
\frac{\zeta_{\Bbb F}(2-s_2)  \zeta_{\Bbb F}(1+s_1-s_2)    }
{\zeta_{\Bbb F}(s_2)   \zeta_{\Bbb F}(2-s_1+s_2)}
\frac{X^{s_1+s_2}Y^{2-s_2}}{s_1s_2(2-s_2)}.
\end{eqnarray*}

We consider the rectangle domain of $s_2$ formed by the four points
  $b_2\pm iT_2$ and $7/4\pm iT_2.$
By the residue theorem we get
\begin{eqnarray}
  \mathfrak{L}_{3}(X,Y,T)= -\mathfrak{L}_{31}(X,Y,T)
+\mathfrak{L}_{32}(X,Y,T)+\mathfrak{L}_{33}(X,Y,T)
\end{eqnarray}
where
\begin{eqnarray*}
&&\mathfrak{L}_{31}(X,Y,T):= \frac{1}{ (2\pi i)^2 }\int_{b_1-iT_1}^{b_1+iT_1}ds_1\int_{b_2+iT_2}^{7/4+iT_2}
\rho_{\Bbb F}
\frac{\zeta_{\Bbb F}(2-s_2)  \zeta_{\Bbb F}(1+s_1-s_2)    }
{\zeta_{\Bbb F}(s_2)   \zeta_{\Bbb F}(2-s_1+s_2)}
\frac{X^{s_1+s_2}Y^{2-s_2}}{s_1s_2(2-s_2)}   ds_2,\\
&&\mathfrak{L}_{32}(X,Y,T):=  \frac{1}{ (2\pi i)^2 }\int_{b_1-iT_1}^{b_1+iT_1}ds_1\int_{7/4-iT_2}^{7/4+iT_2}
\rho_{\Bbb F}
\frac{\zeta_{\Bbb F}(2-s_2)  \zeta_{\Bbb F}(1+s_1-s_2)    }
{\zeta_{\Bbb F}(s_2)   \zeta_{\Bbb F}(2-s_1+s_2)}
\frac{X^{s_1+s_2}Y^{2-s_2}}{s_1s_2(2-s_2)}  ds_2,\\
&&\mathfrak{L}_{33}(X,Y,T):= \frac{1}{ (2\pi i)^2 }\int_{b_1-iT_1}^{b_1+iT_1}ds_1\int_{b_2-iT_2}^{7/4-iT_2}
\rho_{\Bbb F}
\frac{\zeta_{\Bbb F}(2-s_2)  \zeta_{\Bbb F}(1+s_1-s_2)    }
{\zeta_{\Bbb F}(s_2)   \zeta_{\Bbb F}(2-s_1+s_2)}
\frac{X^{s_1+s_2}Y^{2-s_2}}{s_1s_2(2-s_2)}  ds_2.
\end{eqnarray*}

From Lemma 2.2 it is easy to see that
\begin{eqnarray}
\mathfrak{L}_{31}(X,Y,T)\ll \frac{X^2Y}{T^2}\log^3 T+\frac{X^{\frac{11}{4}}Y^{\frac 14}}{T^{\frac 12}}\log^3 T
\end{eqnarray}
and
\begin{eqnarray}
\mathfrak{L}_{33}(X,Y,T)\ll \frac{X^2Y}{T^2}\log^3 T+\frac{X^{\frac{11}{4}}Y^{\frac 14}}{T^{\frac 12}}\log^3 T.
\end{eqnarray}
By (2.4) of  Lemma 2.2 we  can write
\begin{eqnarray*}
 \mathfrak{L}_{32}(X,Y,T)\ll X^{\frac{11}{4}}Y^{\frac 14}\log T
 \times \mathfrak{L}_{32}^{*}(X,Y,T),
\end{eqnarray*}
where
\begin{eqnarray*}
 \mathfrak{L}_{32}^{*}(X,Y,T):=
 \int_{-T_1}^{ T_1}dt_1\int_{-T_2}^{T_2}
 \frac{|\zeta_{\Bbb F}(\frac 14-it_2)  \zeta_{\Bbb F}(\frac 14+\frac{1}{\log X}+i(t_1-t_2)|    }
 {(|t_1|+1)(|t_2|+1)^2}dt_2.
\end{eqnarray*}
By Lemma 2.1 we have
\begin{eqnarray*}
 \mathfrak{L}_{32}^{*}(X,Y,T)\ll
 \int_{-T_1}^{ T_1}dt_1\int_{-T_2}^{T_2}
 \frac{|\zeta_{\Bbb F}(\frac 34+it_2)  \zeta_{\Bbb F}(\frac 34-\frac{1}{\log X}-i(t_1-t_2)|    }
 {(|t_1|+1)(|t_2|+1)^{\frac{3}{2}}}\times (|t_1-t_2|+1)^{\frac 12}dt_2.
\end{eqnarray*}

Similar to (5.11) and (5.12) we have
\begin{eqnarray*}
 \mathfrak{L}_{32}^{*}(X,Y,T)\ll \log^2 T.
\end{eqnarray*}
Thus
\begin{eqnarray}
 \mathfrak{L}_{32}(X,Y,T)\ll X^{\frac{11}{4}}Y^{\frac 14}\log^3 T.
\end{eqnarray}

From (6.22)-(6.25) we get
\begin{eqnarray}
 \mathfrak{L}_{3}(X,Y,T)\ll X^{\frac{11}{4}}Y^{\frac 14}\log^3T
 +\frac{X^2Y}{T^2}\log^3 T.
\end{eqnarray}

\subsection{\bf Upper bound of $\mathfrak{L}_{4}(X,Y,T) $ }\
It is easy to see that
\begin{eqnarray*}
Res_{w=3-s_1-s_2}\mathcal{G}(w;s_1,s_2)&& =\rho_{\Bbb F} \frac{    \zeta_{\Bbb F}(3-s_1-s_2)  \zeta_{\Bbb F}(2-s_2)  \zeta_{\Bbb F}(2-s_1)    }
{\zeta_{\Bbb F}(s_1)\zeta_{\Bbb F}(s_2)  \zeta_{\Bbb F}(4-s_1-s_2)}
\frac{X^{s_1+s_2}Y^{3-s_1-s_2}}{s_1s_2(3-s_1-s_2)}\nonumber\\
&&=m(s_1,s_2),
\end{eqnarray*}
say. By the residue theorem we can write
\begin{eqnarray}
\mathfrak{L}_{4}(X,Y,T)=-\mathfrak{L}_{41}(X,Y,T)+
\mathfrak{L}_{42}(X,Y,T)+\mathfrak{L}_{43}(X,Y,T),
\end{eqnarray}
where
\begin{eqnarray*}
&&\mathfrak{L}_{41}(X,Y,T)=\frac{1}{ (2\pi i)^2 }\int_{b_1-iT_1}^{b_1+iT_1}ds_1\int_{b_2+iT_2}^{\frac 32+iT_2}
m(s_1,s_2)ds_2,\\
&&\mathfrak{L}_{42}(X,Y,T)=\frac{1}{ (2\pi i)^2 }\int_{b_1-iT_1}^{b_1+iT_1}ds_1\int_{\frac 32-iT_2}^{\frac 32+iT_2}
m(s_1,s_2)ds_2,\\
&&\mathfrak{L}_{43}(X,Y,T)=\frac{1}{ (2\pi i)^2 }\int_{b_1-iT_1}^{b_1+iT_1}ds_1\int_{b_2-iT_2}^{\frac 32+iT_2}
m(s_1,s_2)ds_2.
\end{eqnarray*}

Suppose $b_1\leq \sigma_1\leq 3/2, |t_1|\leq T$ and $b_2\leq \sigma_2\leq 3/2, |t_2|\leq 2T.$ By Lemma 2.2 it is easy to check that
\begin{eqnarray}
&&m(s_1,s_2)\ll
 \frac{X^{ \sigma_1+\sigma_2}Y^{3-\sigma_1-\sigma_2}
 \log^2(|t_1|+1)\log^2(|t_2|+1)\log(|t_1+t_2|+1)  }
 {(|t_1|+1)^{2-\sigma_1}(|t_2|+1)^{2-\sigma_2}(|t_1+t_2|+1)^{3-\sigma_1-\sigma_2} }.
\end{eqnarray}
From (6.28) we get ($t_2=T_2=2T$)
\begin{eqnarray}
&&\mathfrak{L}_{41}(X,Y,T) \ll \frac{X^2Y\log^5 T}{T^2}+
\frac{X^{\frac 52}Y^{\frac 12}\log^5 T}{T}
\end{eqnarray}
and ($t_2=-T_2=-2T$)
\begin{eqnarray}
&&\mathfrak{L}_{43}(X,Y,T) \ll \frac{X^2Y\log^5 T}{T^2}+
\frac{X^{\frac 52}Y^{\frac 12}\log^5 T}{T}.
\end{eqnarray}

Now we consider $\mathfrak{L}_{42}(X,Y,T).$ Change the order of
$s_1$ and $s_2$ and then use the residue theorem to $s_1$ we get
\begin{eqnarray}
\mathfrak{L}_{42}(X,Y,T)=-\mathfrak{L}_{421}(X,Y,T)+
\mathfrak{L}_{422}(X,Y,T)+\mathfrak{L}_{423}(X,Y,T),
\end{eqnarray}
where
\begin{eqnarray*}
&&\mathfrak{L}_{421}(X,Y,T)=\frac{1}{ (2\pi i)^2 }\int_{\frac 32-iT_2}^{\frac 32+iT_2}ds_2\int_{b_1+iT_1}^{\frac 32+iT_1}m(s_1,s_2)ds_1,\\
&&\mathfrak{L}_{422}(X,Y,T)=\frac{1}{ (2\pi i)^2 }\int_{\frac 32-iT_2}^{\frac 32+iT_2}ds_2\int_{\frac 32-iT_1}^{\frac 32+iT_1}
m(s_1,s_2)ds_2,\\
&&\mathfrak{L}_{423}(X,Y,T)=\frac{1}{ (2\pi i)^2 }\int_{\frac 32-iT_2}^{\frac 32+iT_2}ds_2\int_{b_1-iT_1}^{\frac 32-iT_1}
m(s_1,s_2)ds_1.
\end{eqnarray*}

From (6.28) we get
\begin{eqnarray}
\ \ \ \ \mathfrak{L}_{421}(X,Y,T)&&\ll
{ }\int_{-T_2}^{ T_2}\frac{1}{(|t_2|+1)^{\frac 12}}dt_2
\int_{b_1 }^{\frac 32 }\frac{X^{\sigma_1+\frac 32}Y^{\frac 32-\sigma_1}\log^5 T}
{T^{2-\sigma_1} (|T+t_2|+1)^{\frac 32-\sigma_1} } d\sigma_1\\
&&\ll \frac{X^{\frac 52}Y^{\frac 12}\log^5 T}{T}
\int_{-T_2}^{T_2}\frac{1}{(|t_2|+1)^{\frac 12}(|T+t_2|+1)^{\frac 12}}dt_2\nonumber\\
&&\ \ \ +\frac{X^3 \log^5 T}{T^{\frac 12} }
\int_{-T_2}^{T_2}\frac{1}{(|t_2|+1)^{\frac 12}   }dt_2\nonumber\\
&&\ll \frac{X^{\frac 52}Y^{\frac 12}\log^5 T}{T }
+ X^3 \log^5 T   \nonumber
\end{eqnarray}
by noting that
\begin{eqnarray*}
&&\int_{-T_2}^{T_2}\frac{1}{(|t_2|+1)^{\frac 12}(|T+t_2|+1)^{\frac 12}}
dt_2\ll 1.
\end{eqnarray*}

Similarly we have
\begin{eqnarray}
 \mathfrak{L}_{423}(X,Y,T)\ll  \frac{X^{\frac 52}Y^{\frac 12}\log^5 T}{T }
+ X^3 \log^5 T.
\end{eqnarray}

Finally we consider $\mathfrak{L}_{422}(X,Y,T).$ Suppose
$s_1=3/2+it_1, s_2=3/2+it_2.$   By Lemma 2.2 we have
\begin{eqnarray*}
m(s_1,s_2)\ll X^3\log^3 T \times \frac{|\zeta_{\Bbb F}(\frac 12+it_1)\zeta_{\Bbb F}(\frac 12+it_2)|}{(|t_1|+1)(|t_2|+1)},
\end{eqnarray*}
which combining with (5.10) implies that
\begin{eqnarray}
\ \ \ \  \mathfrak{L}_{422}(X,Y,T) \ll X^3  \log^9 T.
\end{eqnarray}

From (6.27) and (6.29)-(6.34) we have
\begin{eqnarray}
\mathfrak{L}_{4}(X,Y,T)\ll X^3 \log^9 T
  +\frac{X^{\frac 52}Y^{\frac 12}\log^5 T}{T  }
+\frac{X^2Y\log^5 T}{T^2 }.
\end{eqnarray}

\subsection{\bf Proof of Theorem 2: completion}\

Choose $T=Y^2.$
From (6.1)-(6.5), (6.9), (6.15), (6.21), (6.26), (6.35)   we get
\begin{eqnarray}
\ \ \ \ \ \ C_{\Bbb F,2}(X,Y)&&=\frac{\rho^2_{\Bbb F}}{2\zeta_{\Bbb F}(2)}X^2Y+\frac{\zeta_{\Bbb B}(0)\rho^2_{\Bbb F}}{4\zeta^2_{\Bbb F}(2)}X^4+O\left(\frac{X^2Y^{1+\varepsilon}}{T}+\frac{X^{\frac 52}Y^{\frac 12}\log^7 T}{T}
\right)\\
&&\ \ \ +O\left(X^3\log^9 T+ X^{\frac{11}{4}}Y^{\frac 14}\log^3 T+X^{\frac{3}{2}}Y \log^3 T
\right)\nonumber\\
&& \ \ \ +O\left( X^{\frac{24}{5}}Y^{-\frac 25}+X^{2}Y^{\frac 23} \log^5 T\right)\nonumber\\
 &&=\frac{\rho^2_{\Bbb F}}{2\zeta_{\Bbb F}(2)}X^2Y+\frac{\zeta_{\Bbb B}(0)\rho^2_{\Bbb F}}{4\zeta^2_{\Bbb F}(2)}X^4\nonumber\\
 &&\ \ \ +O\left(X^3\log^9 Y+ X^{\frac{11}{4}}Y^{\frac 14}\log^3 Y+X^{\frac{3}{2}}Y \log^3 Y
\right)\nonumber\\
&& \ \ \ +O\left( X^{\frac{24}{5}}Y^{-\frac 25}+X^{2}Y^{\frac 23} \log^5 Y\right) \nonumber\\
&&=\frac{\rho^2_{\Bbb F}}{2\zeta_{\Bbb F}(2)}X^2Y+\frac{\zeta_{\Bbb B}(0)\rho^2_{\Bbb F}}{4\zeta^2_{\Bbb F}(2)}X^4\nonumber\\
 &&\ \ \ +O\left(X^{\frac{24}{5}}Y^{-\frac 25}+X^{2}Y^{\frac 23} \log^5 Y +X^{\frac{3}{2}}Y \log^3 Y
\right)\nonumber
  \end{eqnarray}
by noting that $Y\geq X^2.$  This completes the proof of Theorem 2.

\begin{flushleft}

\bigskip

Wenguang Zhai\\
Department of Mathematics, \\
China University of Mining and Technology, \\
Beijing 100083, P. R. China\\
e-mail: zhaiwg@hotmail.com
\end{flushleft}


\begin{thebibliography}{99}

\bibitem{CK} T. H. Chan and A.V. Kumchev, On sums of Ramanujan sums, Acta arithm.
{\bf 152} (2012), 1-10.



\bibitem{Gr}A. Grytczuk, On Ramanujan sums on arithmetical semigroups, Tsukuba J.
Math. {\bf 16} (1992), 315-319.

\bibitem{He}E. Hecke, Vorlesungen \"{u}ber die Theorie der algebraischen Zahlen, 2nd ed.,
Chelsea Publ. Co., (New York, 1948).


\bibitem{Iv}A. Ivi\'c, The Riemann Zeta-Function. Theory and Applications, Wiley, New York, 1985.


\bibitem{IK}H. Iwaniec and E. Kowalski, Analytic Number Theory, AMS, 2004.


\bibitem{Kn}J. Knopfmacher, Abstract analytic number theory, North Holland Publ. Co.,
(Amsterdam-Oxford, 1975).


\bibitem{Kr}E. Kr\"{a}tzel, Zahlentheorie, VEB Deutscher Verlag der Wissenschaften, (Berlin,
1981).

\bibitem{Lan}E. Landau,   Einf\"{u}hrung in die elementare und analytische Theorie der algebraischen Zahlen und
der Ideale. 2nd ed. New York 1949.



\bibitem{Mu} W. M\"{u}ller,  On the distribution of ideals in cubic number fields.
 Monatsh. Math. {\bf 106} (1988), 211-219.

\bibitem{No3}W. G. Nowak, On the Distribution of Integer Ideals in Algebraic Number Fields,
Math. Nachr. {\bf 161}(1993), 59-74.


\bibitem{No1}W. G. Nowak, The average size of Ramanujan sums over quadratic number
fields, Arch. Math. {\bf 99} (2012), 433-442.

\bibitem{No2}W. G. Nowak, On Ramanujan sums over the Gaussian integers. Math. Slovaca,
{\bf 63} (2013), No. 4, 725-732.



\bibitem{Te}G. Tenenbaum, {\it Introduction to analytic and probabilistic number theory},
 Third edition, Graduate Studies in Mathematics {\bf 163}, AMS 2015.

\bibitem{Ti}E. C. Titchmarsh, The theory of the Riemann Zeta-Function, Oxford University Press, Oxford, 1951.

\bibitem{TZ}L. T\'{o}th and W. Zhai,  On the average number of the  cyclic subgroups of the groups $\Z_{n_1} \times \Z_{n_2}\times \Z_{n_3}$ with $n_1,n_2,n_3\leq x,$
    Res. Number Theory(2020), 6:12.

\bibitem{Zh}W. Zhai, The average size of Ramanujan sums over quadratic number fields,
Ramanujan J., to appear.



\end{thebibliography}
\end{document}